\documentclass[11pt, eqno]{article}
\usepackage{amssymb,amsmath,latexsym}
\usepackage{epsfig}
\usepackage{color}

\newtheorem{thm}{Theorem}[section]
\newtheorem{lemma}[thm]{Lemma}

\newtheorem{defn}[thm]{Definition}
\newtheorem{prop}[thm]{Proposition}

\newtheorem{remark}[thm]{Remark}

\numberwithin{equation}{section}

\def\pf{{\medskip\noindent {\bf Proof. }}}
\def\qed{{\hfill $\Box$ \bigskip}}
\def\sA {{\mathcal A}}

  \def\sL {{\mathcal L}}

  \def\bR {{\mathbb R}}

\def\R {{\mathbb R}}
\def\RR {{\mathbb R}}
\def\nn{\nonumber}

\def\wt{\widetilde}
\def\wh{\widehat}

\def\E{{\mathbb E}}
\def\P{{\mathbb P}}

\def\bea{\begin{align*}}
\def\eea{\end{align*}}
\def\bee{\begin{equation}}
\def\eee{\end{equation}}

\def\eps{\varepsilon}

\def\wh{\widehat}

\def\1{{\bf 1}}

\makeatletter
\@addtoreset{equation}{section}

\makeatother

\begin{document}
\bibliographystyle{plain}

\title{\Large \bf
Potential theory of subordinate Brownian motions with
Gaussian components
}

\author{{\bf Panki Kim}\thanks{
This research was supported by Basic Science Research Program through the National Research Foundation of Korea(NRF) funded by the Ministry of Education, Science and Technology(0409-20110087) }, \quad {\bf Renming Song}\thanks{Research supported in part by a grant from the Simons
Foundation (208236). }
\quad and \quad {\bf
Zoran Vondra\v{c}ek}\thanks{Research partially supported by the MZOS
grant 037-0372790-2801 of the Republic of Croatia.}}

\date{(April 2, 2012)}

\maketitle

\begin{abstract}
In this paper we study a subordinate Brownian motion with a Gaussian component
and a rather general discontinuous part. The assumption on the subordinator
is that its Laplace exponent is a complete Bernstein function with a L\'evy density satisfying a certain
growth condition near zero. The main result is a boundary Harnack principle with explicit boundary decay rate
for non-negative harmonic functions of the process in $C^{1,1}$ open sets.
As a consequence of the boundary Harnack principle, we establish sharp
two-sided estimates on the Green function of the subordinate Brownian motion in any bounded $C^{1,1}$ open set $D$
and identify the Martin boundary of $D$ with respect to the subordinate Brownian motion with the Euclidean boundary.
\end{abstract}



\bigskip
\section{Introduction}

The infinitesimal generator of a $d$-dimensional rotationally invariant L\'evy process is
a non-local operator of the form $\sL = b \Delta + \sA$ where $b\ge 0$ and
$$
\sA f(x)=\int_{\R^d}\left(f(x+y)-f(x)-\nabla f(x) \cdot y {\mathbf 1}_{\{|y|\le 1\}}\right)\, \nu(dy)
=\lim_{\epsilon\to 0} \int_{\{|y|>\epsilon \}}\left(f(x+y)-f(x)\right)\, \nu(dy)\, .
$$
The measure $\nu$  on $\R^d\setminus \{0\}$ is invariant under rotations around origin
and satisfies $\int_{\R^d}(1\wedge |y|^2)\, \nu(dy) <\infty$. When $\nu=0$, the
operator $\sL$ is proportional to the Laplacian, hence a local operator, while
when $b=0$, the operator $\sL$ is a purely non-local integro-differential operator.
In particular, if $b=0$ and $\nu(dx)=c|x|^{-d-\alpha}dx$, $\alpha \in (0, 2)$,
then $\sA$ is proportional to
the fractional Laplacian
$\Delta^{\alpha/2}:=-(-\Delta)^{\alpha/2}$.
L\'evy processes are of intrinsic importance in probability theory, while integro-differential
operators
 are important in the theory of partial differential equations.
Most of the research in the potential theory of L\'evy processes in
the last fifteen years concentrates on purely discontinuous L\'evy processes,
such as rotationally invariant stable processes,
or equivalently, on purely non-local operators of the type $\sA$.
For summary of some recent results from a probabilistic point of view one can consult
\cite{BBKRSV, C0, KSV2, KSV3}
and references therein. We refer the readers to \cite{CSS, CaS, CV} for a
sample of recent progress in the PDE literature, mostly for the case of a fractional Laplacian
$\Delta^{\alpha/2}$, $\alpha \in (0, 2)$.

In many situations one would like to study operators that have both local and non-local parts.
From a probabilistic point of view, this corresponds to processes with both a Gaussian component
and a jump component. The fact that such a process $X$ has both
Gaussian and jump components is the source of many difficulties in investigating the
potential theory of $X$. The main
difficulty in studying $X$ stems from the fact that it runs on two different scales:
on the small scale the diffusion corresponding to the Gaussian part dominates,
while on the large scale the jumps take over. Another difficulty is encountered when
looking at the exit of $X$ from an open set: for
diffusions, the exit is through the boundary, while for the pure jump processes,
typically the exit happens by jumping out from the
open set. For the process $X$, both cases will occur which makes the process
$X$ much more difficult to study.

Despite the difficulties mentioned above, in the last few years significant progress has been made in
understanding the potential theory of such
processes. Green function estimates (for the whole space) and the Harnack inequality
for a class of processes with both diffusion and
jump components were established in \cite{RSV, SV05}. The parabolic Harnack inequality
and heat kernel estimates were studied in
\cite{SV07} for L\'evy processes in $\bR^d$ that are independent sums of Brownian
motions and symmetric stable processes, and in \cite{CK08} for much more general
symmetric diffusions with jumps. Moreover, an a priori H\"older estimate was established in \cite{CK08}
for bounded parabolic functions. For earlier results on second order integro-differential
operators, one can see \cite{GM} and the references therein.

Important progress has been made in two recent papers \cite{CKSV1, CKSV2}
which consider operators of the type $\Delta + a^{\alpha}\Delta^{\alpha/2}$ for
$a\in [0, M]$. The process corresponding to such an operator is an independent sum of
a Brownian motion and a rotationally invariant $\alpha$-stable process with weight $a$.
In \cite{CKSV1} the authors established a (uniform in $a$) boundary Harnack principle
(BHP) with explicit boundary decay rate for non-negative harmonic functions with respect to
$\Delta + a^{\alpha}\Delta^{\alpha/2}$ in $C^{1,1}$ open sets. By using the BHP, the second
paper \cite{CKSV2} established sharp Green function estimates in bounded $C^{1,1}$ open
sets $D$, and identified the Martin boundary of $D$ for the operator $\Delta +
a^{\alpha}\Delta^{\alpha/2}$ with its Euclidean boundary.

The purpose of the current paper is to extend the results in \cite{CKSV1, CKSV2} to
more general operators than $\Delta + a^{\alpha}\Delta^{\alpha/2}$. Analytically,
the operators that we consider are certain functions of the Laplacian.
To be more precise,
we consider a Bernstein function $\phi:(0,\infty)\to (0,\infty)$ with $\phi(0+)=0$, i.e.,
 $\phi$ is of the form
\begin{equation}\label{e:bernstein-function}
\phi(\lambda)=b \lambda +\int_{(0,\infty)}(1-e^{-\lambda t})\, \mu(dt)\, ,\quad \lambda >0\, ,
\end{equation}
where $b\ge 0$ and $\mu$ is a measure on $(0,\infty)$ satisfying
$\int_{(0,\infty)}(1\wedge t)\, \mu(dt)<\infty$. $\mu$ is called the L\'evy measure of $\phi$.
By Bochner's functional calculus one can define the operator
$\phi(\Delta):=-\phi(-\Delta)$ which on
$C_b^2(\R^d)$, the collection of bounded $C^2$ functions in $\RR^d$ with bounded derivatives,
turns out to be an integro-differential operator of the type
$$
b \Delta f(x)+\int_{\R^d}\left(f(x+y)-f(x)-\nabla f(x) \cdot y {\mathbf 1}_{\{|y|\le 1\}}\right)\, \nu(dy)\, ,
$$
where the measure $\nu$ has the form $\nu(dy)=j(|y|)\, dy $ with $j:(0,\infty)\to (0,\infty)$ given by
$$
j(r)=\int_0^{\infty} (4\pi t)^{-d/2} e^{-r^2/(4t)}\, \mu(dt)\, .
$$
In order for the operator to have both local and non-local parts we will assume that $b>0$
and $\mu\neq 0$. In fact, without loss of generality, throughout the paper we always
suppose that  $b=1$. Note that by taking $\phi(\lambda)=\lambda + a^{\alpha}
\lambda^{\alpha/2}$ we are back to the operator $\Delta + a^{\alpha}\Delta^{\alpha/2}$.

The operator $\phi(\Delta)$ is the infinitesimal generator of the L\'evy process
$X$ that can be constructed as follows: Recall that a one-dimensional  L\'evy
process $S=(S_t: t\ge 0)$ is called a subordinator if
it is non-negative and $S_0=0$.
A subordinator $S$ can be characterized by its Laplace exponent $\phi$ through the equality
$$
\E[e^{-\lambda S_t}]=e^{-t \phi(\lambda)}, \quad t>0, \lambda>0\, .
$$
The Laplace exponent $\phi$ can be written in the form \eqref{e:bernstein-function}.
We will assume that $b=1$.
Suppose that $W=(W_t: t\ge 0)$ is a $d$-dimensional
Brownian motion and $S=(S_t: t\ge 0)$ is a subordinator, independent of $W$, with Laplace exponent $\phi$.
The process $X=(X_t:\, t\ge0)$ defined by $X_t=W_{S_t}$ is
called a subordinate Brownian motion and its infinitesimal generator is
$\phi(\Delta)$. It is a sum of a Brownian motion and an independent purely
discontinuous (rotationally invariant) L\'evy process.

Potential theory of one-dimensional subordinate Brownian motions in this setting
was studied in \cite{KSV1}. In the current paper we look at the case when
$d\ge 2$. In order for our methods to work we need additional assumptions on the
Bernstein function $\phi$. We will assume that $\phi$ is a complete Bernstein
function, namely that the L\'evy measure $\mu$ has a completely monotone density.
By a slight abuse of notation we will denote the density by $\mu(t)$. For the
L\'evy density $\mu$ we assume a growth condition near zero: For any $K>0$, there exists $c=c(K) >1$ such that
\begin{equation}\label{H:1a}
\mu(r)\le c\, \mu(2r), \qquad \forall r\in (0, K)\, .
\end{equation}
We will later explain the role of these additional assumptions.

To state our main result, we first recall that an open set $D$ in $\bR^d$
(when $d\ge 2$) is said to be $C^{1,1}$ if
there exist a localization radius $R>0$ and a constant $\Lambda >0$ such that
for every $Q\in \partial D$, there exist a
$C^{1,1}$-function $\varphi=\varphi_Q: \bR^{d-1}\to \bR$ satisfying $\varphi (0)=0$,
$ \nabla\varphi (0)=(0, \dots 0)$, $\| \nabla
\varphi  \|_\infty \leq \Lambda$, $| \nabla \varphi (x)-\nabla \varphi (y)| \leq
\Lambda |x-y|$, and an orthonormal coordinate system $CS_Q$: $y=(y_1, \cdots,
y_{d-1}, y_d)=:(\wt y, \, y_d)$ with its origin at $Q$ such that
$$
B(Q, R)\cap D=\{ y=(\wt y, y_d)\in B(0, R) \mbox{ in } CS_Q: y_d >
\varphi (\wt y) \}.
$$
The pair $(R, \Lambda)$ is called the characteristics of the $C^{1,1}$ open set $D$.
Note that a $C^{1,1}$ open set can be unbounded and disconnected.

For any $x\in D$, $\delta_D(x)$ denotes the Euclidean distance
between $x$ and $D^c$. For any $x\notin D$, $\delta_{\partial D}(x)$
denotes the Euclidean distance between $x$ and $\partial D$.
It is well known that, if $D$ is a $C^{1, 1}$ open set $D$ with
characteristics $(R, \Lambda)$, there exists $\wt R \le  R$ such
that $D$ satisfies both the {\it uniform interior ball condition} and the
{\it uniform exterior ball condition} with radius $\wt R$:
for every $x\in D$ with $\delta_{D}(x)<\wt R$ and $y\in \bR^d
\setminus \overline D$ with $\delta_{\partial D}(y)<\wt R$, there are $z_x,
z_y\in \partial D$ so that $|x-z_x|=\delta_{ D}(x)$,
$|y-z_y|=\delta_{\partial D}(y)$ and that $B(x_0,\wt R)\subset D$ and $B(y_0,
\wt R)\subset \bR^d \setminus \overline D$ where $x_0=z_x+\wt
R(x-z_x)/|x-z_x|$ and $y_0=z_y+\wt R(y-z_y)/|y-z_y|$. Without loss
of generality, throughout this paper, we assume that the
characteristics $(R, \Lambda)$ of a $C^{1, 1}$ open set satisfies
$R=\wt R\le 1$ and $\Lambda\ge 1$.

For any open set $D\subset \bR^d$, $\tau_D:=\inf\{t>0: \,
X_t\notin D\}$ denotes the first exit time from $D$ by $X$.

\begin{defn}\label{D:1.1}
A function $f: \R^d\mapsto [0, \infty)$ is said to be
\begin{description}
\item{(1)} harmonic in an open set $D\subset \R^d$ with respect to
$X$ if for every open set $B$ whose closure is a compact subset of $D$,
\begin{equation}\label{e:har}
f(x)= \E_x \left[ f(X_{\tau_{B}})\right] \qquad
\mbox{for every } x\in B;
\end{equation}
\item{(2)} regular harmonic in $D$ for $X$ if
for each $x \in D$,
$f(x)= \E_x\left[f(X_{\tau_{D}})\right]$.
\end{description}
\end{defn}
We note that, by the strong Markov property of $X$,
every regular harmonic function is automatically harmonic.

Let $Q\in \partial D$. We will say that a function $f:\bR^d\to \R$
vanishes continuously on $ D^c \cap B(Q, r)$ if $f=0$ on $ D^c \cap
B(Q, r)$ and $f$ is continuous at every point of $\partial D\cap
B(Q,r)$. The following is the main result of this paper.

\begin{thm}\label{t:main}
Suppose that the Laplace exponent $\phi$ of
the subordinator $S$, independent of the Brownian motion $W$,  is a complete Bernstein
function and that the L\'evy density of $S$ satisfies
\eqref{H:1a}. Let $X=(X_t:\, t\ge0 )$ be the subordinate Brownian motion defined by $X_t=W(S_t)$. For any $C^{1, 1}$ open set $D$ in $\bR^d$ with characteristics
$(R, \Lambda)$, there exists a positive constant
$C=C(d, \Lambda, R)$ such that for  $r \in (0, R]$, $Q\in \partial D$ and any
nonnegative function $f$ in $\R^d$ which is
harmonic in $D \cap B(Q, r)$ with respect to $X$
 and vanishes
continuously on $ D^c \cap B(Q, r)$, we have
\begin{equation}\label{e:bhp_m}
\frac{f(x)}{\delta_D(x)}\,\le C\,\frac{f(y)}{\delta_D(y)} \qquad
\hbox{for every } x, y\in  D \cap B(Q, r/2).
\end{equation}
\end{thm}

We note that \eqref{e:bhp_m} is a strengthened version of the usual boundary
Harnack principle  stated for the ratio of two non-negative functions,
$f$ and $g$, harmonic in $D \cap B(Q, r)$ with respect to $X$, and which says that
$$
\frac{f(x)}{g(x)}\,\le C\,\frac{f(y)}{g(y)} \qquad
\hbox{for every } x, y\in  D \cap B(Q, r/2).
$$
Indeed, the above inequality is a consequence of \eqref{e:bhp_m}.
We note  that \eqref{e:bhp_m} gives the precise boundary decay of non-negative harmonic functions and
 that the function $x\mapsto \delta_D(x)$ is not harmonic in $D \cap B(Q, r)$ with respect to $X$.

\begin{remark}\label{r:counterexample}{\rm
The same type of boundary Harnack principle in $C^{1,1}$ domains is
valid also for Brownian motions, namely the boundary decay rate is of
the order $\delta_D(x)$. Since on the small scale the diffusion part of
$X$ dominates, one would expect that
harmonic functions of $X$ and of Brownian
motion
 have the same decay rate
at the boundary.
For this reason, some people might expect that some kind of perturbation methods
can be used to prove the BHP for $X$.
We note that it is unlikely that any perturbation
method would work because of the following: Suppose that instead of $X$ we consider a process
$X^a$ with the infinitesimal generator
$$
{\sL}^af(x)= \Delta f(x) +\int_{\R^d}\left(f(x+y)-f(x)-\nabla f(x) \cdot y
{\mathbf 1}_{\{|y|\le 1\}}\right)\, {\nu}
^a(dy)\, ,
$$
where ${\nu}^a(dy)={\mathbf 1}_{\{|y|\le
a\}}\, \nu(dy)$ with $0<a<\infty$.
Thus  $X^a$ is the process $X$ with jumps of size larger than
$a$ suppressed.
 In Section \ref{counterexample} we present an example of a
(bounded) $C^{1,1}$ domain $D$ on which the boundary Harnack principle for
$X^a$ fails, even for regular harmonic functions vanishing on $D^c$. Note that
if we think of $X$ as a perturbation of Brownian motion, then
$X^a$ is an even smaller perturbation of the same Brownian motion. The counterexample
in Section \ref{counterexample} shows that, despite the (seemingly)
local nature of the BHP, one needs some information of the structure of large jumps of $X$.
}
\end{remark}

For any open set $D\subset \R^d$, we will use $X^D$ to denote the process
defined by $X^D_t(\omega)=
X_t(\omega)$ if $t<\tau_D(\omega)$ and $X^D_t(\omega)=\partial$ if $t\ge
\tau_D(\omega)$, where $\partial$ is a cemetery point.
The Green function of $X^D$ will be denoted by $G_D(x, y)$. For the precise
definition of $G_D$, see Section 2.

To state our result on Green function estimates, we introduce a function $g_D$ first.
For $d\ge 2$, we define for $x, y \in D$,
$$
g_D(x, y)=
\begin{cases}
 \frac{1} {|x-y|^{d-2}} \left(1\wedge \frac{  \delta_D(x) \delta_D(y)}{ |x-y|^{2}}\right) & \text{ when } d \ge3,\\
 \log\left(1+\frac{  \delta_D(x) \delta_D(y)}{ |x-y|^{2}}\right) & \text{ when }  d=2.
\end{cases}
$$
\begin{thm}\label{t-main-green}
Suppose that the Laplace exponent $\phi$ of $S$ is a complete
Bernstein function and that the L\'evy density of $S$ satisfies
\eqref{H:1a}.
For any bounded $C^{1,1}$ open set $D\subset \R^d$, there exists
$C=C(D)>1$ such that for
all $x, y \in D$
\begin{equation}\label{e:1.3}
C^{-1} \, g_D(x, y) \leq G_D(x, y) \leq C\, g_D(x, y) .
\end{equation}
\end{thm}

Finally, we state the result about the Martin boundary of a bounded $C^{1,1}$ open set $D$ with respect to $X^D$. Fix $x_0\in D$ and define
$$
M_D(x, y):=\frac{G_D(x, y)}{G_D(x_0, y)}, \qquad x, y\in D,~y\neq x_0.
$$
A function $f$ is called a harmonic function for $X^D$ if it is harmonic for $X$ in $D$ and
vanishes outside $D$.
A positive harmonic function $f$ for  $X^{D}$ is minimal if, whenever
$g$ is a positive harmonic function for $X^{D}$ with $g\le f$ on $D$,
one must have $f=cg$ for some constant $c$.
\begin{thm}\label{t-main-martin}
Suppose that $D$ is a bounded $C^{1,1}$ open set in $\R^d$. For every $z\in \partial D$, there exists $M_D(x,z):=\lim_{y\to z}M_D(x,y)$.
Further, for every $z \in \partial D$, $M_D(\cdot, z)$ is a minimal harmonic function for $X^{D}$ and $M_D(\cdot, z_1)\not=M_D(\cdot, z_1)$ if $z_1\not=z_2$.
 Thus the minimal Martin boundary of $D$ can be identified with the Euclidean boundary.
\end{thm}

Thus, by the general theory of Martin representation in \cite{KW} and
Theorem \ref{t-main-martin} we conclude that,  for every  harmonic
function $u\geq 0$ with respect to $X^{D}$, there is a  unique
 finite measure $\nu$ on $\partial
D$ such that
$
 u(x) =\int_{\partial D} M_D(x,z) \nu (dz).
$

Let us now describe the main ingredients of the proof of Theorem \ref{t:main},
the boundary Harnack principle. We follow the general strategy for proving the
boundary Harnack principle in different settings which requires the Carleson estimate, and upper and lower
estimates on exit probabilities and exit times from certain sets usually called
``boxes'' (see \cite{BBC, BB, CKSV1, G, KS}) .
In Theorem \ref{carleson} we prove the Carleson estimate for a Lipschitz
open set by modifying the proof in \cite{CKSV1}. In order to obtain the upper exit
probability and exit times  estimates, we follow the approach from \cite{CKSV1},
the so-called ``test function'' method (which was modeled after some earlier ideas,
see \cite{BBC, G}), but have to make major modifications. In \cite{CKSV1}, the test
functions are power functions of the form $x\mapsto (x_d)^p$ which are either
superharmonic or subharmonic for the corresponding process, and the values of
the generator on these test functions are computed in detail. In our setting,
the power functions are neither superharmonic nor subharmonic, and explicit
calculations cannot be carried out because of the lack of explicit form of the
L\'evy measure. Instead we use the approach developed in \cite{KSV3} for the
case of certain pure-jump subordinate Brownian motions, which seems to be
quite versatile to cover various other cases.

One of the main ingredients
 in \cite{KSV3} comes from the fluctuation theory of one-dimensional
L\'evy processes. Its purpose is to identify a correct boundary decay rate by finding an
appropriate harmonic function. Let $Z=(Z_t:\, t\ge 0)$ be the one-dimensional subordinate
Brownian motion defined by $Z_t:=W^d_{S_t}$, and let $V$ be the renewal function of the
ladder height process of $Z$. The function $V$ is harmonic for the process $Z$ killed upon
exiting $(0,\infty)$, and the function $w(x):= V(x_d){\bf 1}_{\{x_d >0\}}$, $x\in \R^d$, is
harmonic for the process $X$ killed upon exiting the half space $\R^d_+:=\{ x=(x_1, \dots,
x_{d-1}, x_d) \in \bR^d: x_d > 0 \}$ (Theorem \ref{t:Sil}). Therefore, $w$ gives the correct
rate of decay of harmonic functions near the boundary of $\R^d_+$. We will use the function
$w$ as our test function. Note that the assumption that $\phi$ is a complete Bernstein
function implies that $w$ is smooth. Using smoothness and harmonicity of $w$ together
with the characterization of harmonic functions recently established in \cite{C}, we
show that $(\Delta + \mathcal A) w \equiv 0$ on the half space (Theorem \ref{c:Aw=0}).
Consequently we prove the following fact in Lemma \ref{L:Main}, which is the key to
proving upper estimates: If $D$ is a $C^{1,1}$ open set with characteristics $(R, \Lambda)$,
$Q\in \partial D$ and $h(y)=V(\delta_D(y)) \1_{D\cap B(Q,R)}$, then $(\Delta+\mathcal A) h(y)$
is a.e.~well defined and bounded for $y\in D$ close enough to the boundary point $Q$. Using
this lemma, we give necessary exit distribution estimates  in Lemma \ref{L:2}. Here we modify
the test function $h$ by adding a quadratic type function (in one variable) -- this is
necessary due to the presence of the Laplacian. The desired exit distribution estimates
are directly derived by applying Dynkin's formula to the new test function.
The reader will note that our proof is even shorter than the one in \cite{CKSV1},
partly because, in \cite{CKSV1}, the uniformity of the boundary Harnack principle for $\Delta +
a^{\alpha}\Delta^{\alpha/2}$ in the weight $a \in (0, M]$ is established.

In order to prove the lower bound for the exit probabilities  we compare the
process $X$ killed upon exiting a certain box $\wh{D}$ with the so-called
subordinate killed Brownian motion obtained by first killing Brownian motion
upon exiting the box $\wh{D}$, and then by subordinating the obtained process.
If the latter process is denoted by $Y^{\wh{D}}$, then its infinitesimal generator is equal to
$-\phi(-\Delta|_{\wh{D}})$.
 Here $\Delta|_{\wh{D}}$ is the Dirichlet Laplacian and
$-\phi(-\Delta|_{\wh{D}})$ is constructed by Bochner's subordination. The
advantage of this approach is that the exit probabilities of $X^{\wh{D}}$
dominate from the above those of the process $Y^{\wh{D}}$, while the
latter can be rather easily computed, see \cite{SV08}. This idea is
carried out in Lemma \ref{L:200} (as well as for some other lower bounds
throughout the paper).

Once the boundary Harnack principle has been established, proofs of Theorems
\ref{t-main-green} and \ref{t-main-martin} are similar to the corresponding
proofs in \cite{CKSV2} for the operator $\Delta+ a^{\alpha}
\Delta^{\alpha}$.
Therefore we do not give complete proofs of these two theorems in this paper,
only indicate the necessary changes to the proofs in \cite{CKSV2}.

The rest of the paper is organized as follows.
In the next section we precisely describe
the settings and recall necessary preliminary results. Section 3 is
devoted to the analysis of the process and harmonic functions in the
half-space.
Section 4 is on the analysis in $C^{1,1}$ open sets, and is central to the paper,
and this is where most of the new ideas appear.
In this rather technical section we establish the upper and lower
bounds on the exit probabilities and exit times. In Section 5 we
first prove the Carleson estimate for Lipschitz open sets and then
prove the main Theorem \ref{t:main}. In Section 6 we provide the
counterexample already mentioned in Remark \ref{counterexample}.
Finally, in Section 7 we explain the differences between the proofs of
Theorems \ref{t-main-green} and \ref{t-main-martin} and their
counterparts from \cite{CKSV2}.

Throughout this paper, the constants $C_1$, $C_2$,  $R$, $R_1$, $R_2$, $R_3$ will be fixed.
The lowercase constants $c_1, c_2, \cdots$ will denote generic constants whose exact values
are not important and can change from one appearance to another.
The dependence of the lower case constants on the  dimension $d$
may not be mentioned explicitly. We will
use ``$:=$" to denote a definition, which is read as ``is defined to
be". For $a, b\in \bR$, $a\wedge b:=\min \{a, b\}$ and $a\vee
b:=\max\{a, b\}$. For every function $f$, let $f^+:=f \vee 0$.
For every
function $f$, we extend its definition to the cemetery point $\partial$ by setting
$f(\partial )=0$.
We will use $dx$ to denote the
Lebesgue measure in $\bR^d$ and, for a Borel set $A\subset \bR^d$, we
also use $|A|$ to denote its Lebesgue measure.

\section{Setting and Preliminary Results}

A $C^{\infty}$ function $\phi:(0,\infty)\to [0,\infty)$ is called a
Bernstein function if $(-1)^n D^n \phi\le 0$ for every
$n=1, 2, \dots$.
Every Bernstein function has a representation $
\phi(\lambda)=a+b\lambda +\int_{(0,\infty)}(1-e^{-\lambda t})\,
\mu(dt) $ where $a,b\ge 0$ and $\mu$ is a measure on $(0,\infty)$
satisfying $\int_{(0,\infty)}(1\wedge t)\, \mu(dt)<\infty$; $a$ is
called the killing coefficient, $b$ the drift and $\mu$ the L\'evy
measure of the Bernstein function. A Bernstein function $\phi$ is
called a complete Bernstein function if the L\'evy measure $\mu$ has
a completely monotone density $\mu(t)$, i.e., $(-1)^n D^n \mu(t)\ge 0$
for every non-negative integer $n$ and all $t>0$. Here and below, by abuse of
notation  we denote the L\'evy density by $\mu(t)$. For more on
Bernstein and complete Bernstein functions we refer the readers to
\cite{SSV}.

A Bernstein function $\phi$ on $(0, \infty)$ is the Laplace exponent of a
subordinator if and only if $\phi(0+)=0$. Suppose that $S$ is a
subordinator with Laplace exponent $\phi$. $S$ is called a complete
subordinator if $\phi$ is a complete Bernstein function.
The potential measure $U$ of $S$ is defined by
\begin{equation}\label{potential measure}
U(A)=\E \int_0^{\infty}
{\bf 1}_{\{S_t\in A\}}
\, dt, \quad A\subset [0,
\infty).
\end{equation}
Note that $U(A)$ is the expected time the subordinator $S$ spends in
the set $A$.

Throughout the remainder of this paper, we assume that  $S=(S_t:\ t\ge 0)$
is a complete subordinator with a positive drift and, without loss
of generality, we shall assume that the drift of $S$ is equal to 1.
Thus the Laplace exponent of $S$ can be written as
$$
\phi(\lambda):=\lambda+\psi(\lambda)
\qquad
\text{ where }
~~
\psi(\lambda):=\int_{(0,\infty)}(1-e^{-\lambda t})\, \mu(dt).
$$
We will exclude the trivial case of $S_t=t$, that is
the case of $\psi\equiv 0$.
Since the drift of $S$ is equal to 1, the potential measure $U$ of $S$ has a
completely monotone density $u$
(cf. \cite[Corollary 5.4 and Corollary 5.5]{BBKRSV}).

Suppose that $W=(W_t: t\ge 0)$ is a Brownian motion in $\R^d$ independent
of $S$ and with
$$
\E_x[e^{i\theta\cdot (W_t-W_0)}]=e^{-t|\theta|^2}, \quad \mbox{ for all }
x, \theta\in \R^d.
$$
The process $X=(X_t: t\ge 0)$ defined by $X_t=W_{S_t}$ is called
a subordinate Brownian motion. It follows from \cite[Chapter 5]{BBKRSV} that
$X$ is a L\'evy process with L\'evy exponent $\phi(|\theta|^2)=|\theta|^2
+\psi(|\theta|^2)$:
$$
\E_x[e^{i\theta\cdot (X_t-X_0)}]=e^{-t\phi(|\theta|^2)}, \quad \mbox{ for all }
x, \theta\in \R^d.
$$
The L\'evy measure of the process $X$ has a density $J$, called
the L\'evy density, given by
 $J(x)=j(|x|)$ where
\begin{equation}\label{e:representation-j}
j(r)
:=\int^{\infty}_0(4\pi t)^{-d/2}e^{-r^2/(4t)}\mu(t)\, dt, \qquad r>0.
\end{equation}
Note that the function $r\mapsto j(r)$ is continuous and decreasing
on $(0, \infty)$. We will sometimes use the notation $J(x,y)$ for
$J(x-y)$.

The function $J(x,y)$ is the L\'evy intensity of $X$. It determines
a L\'evy system for $X$, which describes the jumps of the process
$X$: For any non-negative measurable function $f$ on $\bR_+ \times
\bR^d\times \bR^d$ with $f(s, y, y)=0$ for all $y\in \bR^d$, any
stopping time $T$ (with respect to the filtration of $X$) and any
$x\in \bR^d$,
\begin{equation}\label{e:levy}
\E_x \left[\sum_{s\le T} f(s,X_{s-}, X_s) \right]= \E_x \left[
\int_0^T \left( \int_{\bR^d} f(s,X_s, y) J(X_s,y) dy \right)
ds \right].
\end{equation}
(See, for example, \cite[Proof of Lemma 4.7]{CK1} and \cite[Appendix
A]{CK2}.)

Recall that for any open set $U\subset \bR^d$,
$\tau_U=\inf\{t>0: \, X_t\notin U\}$ is  the first exit time
from $U$ by $X$.
The following simple result will be used in Section 5.

\begin{lemma}\label{L:2.00}
For every $\varrho>0$, there  exists $c=c(\varrho)>0$ such that for
every $x_0 \in \bR^d$ and $r \in (0, \varrho]$,
 \bee \label{e:ext}
c^{-1} r^2 \, \le \, \E_{x_0}\left[\tau_{B(x_0,r)}\right]\, \le\, c\,r^2.
 \eee
\end{lemma}

\pf
In the case $d\ge 3$, this lemma has been proved in \cite{RSV}.
Moreover, one can easily adapt the proofs of \cite[Lemmas 2.1--2.1]{SV05} to arrive at
the desired lower bound for all dimensions.
Here
we provide a proof of the desired upper bound that works for all dimensions.

Let
$C^2_0(\bR^d)$ be the collection of $C^2$ functions in $\RR^d$ vanishing at
infinity.
For any $f\in C^2_0(\R^d)$,
we define
$$
{\sA}f(y)=\int_{\R^d}(f(z+y)-f(y)-z\cdot \nabla f(y)1_{\{|z|<1\}})J(z)dz, \quad y\in \R^d.
$$
Let $N\ge 5$ be such that
$$
\int_{\{|z|<1\}}|z|^2J(z)dz<\frac12N^2.
$$
Let $g$ be a radial $C^2$ function taking values in $[0, 2]$ such that
$$
g(y)=\begin{cases}|y|^2, &|y|<1\\
2, &2\le |y|\le 3\\
0, &|y|>N-1.
\end{cases}
$$
For any $r>0$, put $f(y)=g(y/r)$. Then for $y\in B(0, r)$, $\Delta f(y)=2dr^{-2}$. For any
$y\in B(0, r)$, we have
\begin{eqnarray*}
|{\sA}f(y)|
&=&|\int_{\R^d}(f(z+y)-f(y)-z\cdot \nabla f(y)1_{\{|z|<Nr\}})J(z)dz|\\
&\le&|\int_{\{|z|<Nr\}}(f(z+y)-f(y)-z\cdot \nabla f(y))J(z)dz|
\ + \ \int_{\{Nr<|z|\}}f(y)J(z)dz\\
&\le& c_1r^{-2}\int_{\{|z|<Nr\}}|z|^2J(z)dz  + r^{-2}|y|^2\int_{\{Nr<|z|\}}J(z)dz\\
&\le& c_1r^{-2}\int_{\{|z|<Nr\}}|z|^2J(z)dz +N^{-2}r^{-2}\int_{\{Nr<|z|<1\}}|z|^2J(z)dz+\int_{\{|z|>1\}}J(z)dz.
\end{eqnarray*}
Thus we know that there exist $r_0\in (0, 1)$ and $c_2>0$ such that for any $r\in (0, r_0)$,
$$
\Delta f(y)+{\sA}f(y)\ge c_2 r^{-2}, \quad y\in B(0, r).
$$
Using this and the fact that $\Delta+{\sA}$ is the infinitesimal generator of
the process $X$, by the Dynkin's formula,
we have that for $r\in (0, r_0)$,
\begin{eqnarray*}
\E_{x_0}[\tau_{B(x_0, r)}]&=&\E_{0}[\tau_{B(0, r)}]=\lim_{t\uparrow\infty}\E_{0}[\tau_{B(0, r)}\wedge t]\\
&\le& c_2^{-1} r^2\lim_{t\uparrow\infty}\E_{0}\int^{\tau_{B(0, r)}\wedge t}_0 (\Delta + {\sA})f(X_s)ds\\
&=&c_2^{-1} r^2\lim_{t\uparrow\infty}\E_{0}[f(X_{\tau_{B(0, r)}\wedge t})]\le 2c_2^{-1} r^2\, ,.
\end{eqnarray*}
Now the desired upper bound follows easily.
\qed

In the remainder of this paper, we will need some control on the behavior of $j$ near the origin.
For this, we will assume that for any $K>0$, there exists $c=c(K) >1$ such that
\begin{equation}\label{H:1b}
\mu(r)\le c\, \mu(2r), \qquad \forall r\in (0, K).
\end{equation}
On the other hand, since $\phi$ is a complete Bernstein function, it follows from
\cite[Lemma 2.1]{KSV3} that there exists $c>1$ such that $\mu(t)\le c
\mu(t+1)$ for every  $t>1$. Thus by repeating the proof of
\cite[Lemma 4.2]{RSV} (see also \cite[Proposition 1.3.5]{KSV2}), we can show that
for any $K>0$,
there exists $c=c(K)>1$ such that
\begin{equation}\label{H:1}
j(r)\le c\, j(2r), \qquad \forall r\in (0, K),
\end{equation}
and, there exists $c>1$ such that
\begin{equation}\label{H:2}
j(r)\le c\, j(r+1), \qquad \forall r>1.
\end{equation}
Note that, as a consequence of \eqref{H:1}, we have
that, for any $K>0$,
\begin{equation}\label{H:1n}
j(ar)\le c\, a^{-\nu} j(r), \qquad \forall r\in (0, K)
\quad\text{and}\quad  a \in (0, 1)
\end{equation}
where $c=c(K)$ is the constant in \eqref{H:1} and $\nu=\nu(K):=\log_2 c$.

The following Harnack inequality
will be used to prove the main result of this paper.

\begin{prop}
[Harnack inequality]\label{uhp}
There exists a constant $c>0$ such that for any
$r\in (0,1]$ and $x_0\in \R^d$ and any function $f$ which is nonnegative in
$\R^d$ and harmonic in $B(x_0, r)$ with respect to $X$ we have
$$
f(x)\le c f(y) \qquad \mbox{ for all } x, y\in B(x_0, r/2).
$$
\end{prop}

\pf
We first deal with the case $d\ge 3$.
When $f$ is bounded, this proposition is just \cite[Theorem 4.5]{RSV}.
Using the same argument as in the proof of
\cite[Corollary 4.7]{RSV}, one can easily see that \cite[Theorem
4.5]{RSV} can be extended to  any nonnegative harmonic function.

The assertions of the proposition in the cases of $d=2$ and $d=1$
follow easily from the assertion in the case $d\ge 3$. Since the arguments
are similar, we will only spell out the details in the case $d=2$.
For any $x\in \R^3$, $x=(\wt{x},x^3)$, where $\wt{x}\in \R^2$.
Analog notation will be used also for other objects in $\R^3$.
Let $X=(X_t, \P_x)$ be the subordinate Brownian motion in $\R^3$ and write $X=(\wt{X}, X^3)$.
Note that $\wt{X}$ has the same distribution under
$\P_{(\wt{x},0)}$ and $\P_{(\wt{x},x^3)}$ for any $x^3\in \R$.
Hence we can define $\P_{\wt{x}}:=\P_{(\wt{x},0)}$. The process
$(\wt{X}, \P_{\wt{x}})$ is a subordinate Brownian motion in $\R^2$ via
the same subordinator as the one used to define $X$.
For any given
$\wt{f}:\R^{2}\to [0,\infty)$, we extend it to $\R^3$ by defining
$f(x)=f((\wt{x},x^3)):=\wt{f}(\wt{x})$. Then
\begin{description}
\item{(1)} If $\wt{f}$ is regular harmonic (with respect to $\wt{X}$)
in an open set $\wt{D}\subset \R^{2}$, then $f$ is regular harmonic
(with respect to $X$) in the cylinder $D:=\wt{D}\times \R$.
Indeed, let $\wt{\tau}_{\wt{D}}:=\inf\{t>0:\, \wt{X}_t\notin \wt{D}\}$
be the exit time of $\wt{X}$ from $\wt{D}$, and
$\tau_D:=\inf\{t>0:\, X_t\notin D\}$. Then clearly $\wt{\tau}_{\wt{D}}=\tau_D$. Thus, for any $x=(\wt{x},x^3)\in D$,
$$
\E_x[f(X_{\tau_D})]=\E_{\wt{x}}[\wt{f}(\wt{X}_{\wt{\tau}_{\wt{D}}})]=\wt{f}(\wt{x})=f(x)\, .
$$

\item{(2)} If $\wt{f}$ is harmonic (with respect to $\wt{X}$) in an open set
$\wt{D}\subset \R^{2}$, then $f$ is  harmonic (with respect to $X$) in the
cylinder $D:=\wt{D}\times \R$. Indeed, let $B\subset D$ be open and
relatively compact. Then there exists a cylinder $C=\wt{C}\times \R$
such that $B\subset C$ and $\wt{C}\subset \wt{D}$ is open and relatively
compact (in $\wt{D}$). Since $\wt{f}$ is harmonic (with respect to
$\wt{X}$) in $\wt{D}$, it is regular harmonic in $\wt{C}$. By (1),
$f$ is regular harmonic (with respect to $X$) in $C$, and therefore
also harmonic in $C$. Since $B$ is compactly contained in $C$, we see that
$$
f(x)=\E_x[f(X_{\tau_B})]\, , \qquad \textrm{for all } x\in B\, .
$$
\end{description}

Let $r\in (0,1)$, $\wt{x}_0\in
\R^2$ and define $x_0:=(\wt{x}_0,0)$.
Assume that $\wt{f}:
\R^2\to [0,\infty)$ is harmonic
(with respect to $\wt{X}$) in $B(\wt{x}_0, r)$. Then $f$
defined by $f(x)=\wt{f}(\wt{x})$ is harmonic in $B(\wt{x}_0, r)\times \R$.
In particular, $f$ is harmonic in $B(x_0,r)$. By the assertion in the case $d=3$,
$$
f(x)\le c f(y)\, ,\qquad \textrm{for all }x,y\in B(x_0,r/2)\, .
$$
Let $\wt{x}, \wt{y}\in B(\wt{x}_0,r/2)$, and define $x:=(\wt{x},0)$, $y:=(\wt{y},0)$. Then
$$
\wt{f}(\wt{x})=f(x)\le cf(y)=c\wt{f}(\wt{y})\, .
$$
\qed

It follows from \cite[Chapter 5]{BBKRSV} that the process $X$ has a transition
density $p(t, x, y)$, which is jointly continuous. Using this and the strong Markov property, one can easily check that
$$
p_D(t, x, y):=p(t, x, y)-
\E_x[ p(t-\tau_D, X_{\tau_D}, y); t>\tau_D], \quad x, y \in D
$$
is continuous and is the transition density of $X^D$.
 For any bounded open set $D\subset \R^d$, we
will use $G_D$ to denote the Green function of $X^D$, i.e.,
$$
G_D(x, y):=\int^\infty_0p_D(t, x, y)dt, \quad x, y\in D.
$$
Note that $G_D(x,y)$ is continuous on $\{(x,y) \in D \times D: x\not=y\}$.

\section{Analysis on half-space}

Recall that $X=(X_t:\, t\ge 0)$ is the $d$-dimensional subordinate
Brownian motion defined by $X_t=W_{S_t}$, where $W=(W^1,\dots, W^d)$
is a $d$-dimensional Brownian motion and $S=(S_t:\, t\ge 0)$ an
independent complete subordinator whose drift is equal to 1 and
whose L\'evy density satisfies \eqref{H:1a}.

Let $Z=(Z_t:\, t\ge 0)$ be the one-dimensional subordinate Brownian
motion defined as $Z_t:=W^d_{S_t}$. Let $\overline{Z}_t:=\sup\{0\vee
Z_s:0\le s\le t\}$ be the supremum process of $Z$ and let $L=(L_t:\,
t\ge 0)$ be a local time of $\overline{Z}-Z$ at $0$. $L$ is also
called a local time of the process $Z$ reflected at the supremum.
The right continuous inverse $L^{-1}_t$ of $L$ is a subordinator and
is called the ladder time process of $Z$. The process
$H_t=\overline{Z}_{L^{-1}_t}$ is also a subordinator and is called
the ladder height process of $Z$. (For the basic properties of the
ladder time and ladder height processes, we refer our readers to
\cite[Chapter 6]{Ber}.) The ladder height process $H$ has a drift
(\cite[Lemma 2.1]{KSV1}). The potential measure
of the subordinator $H$ will be denoted by $V$.
Let $V(t):=V((0, t))$ be the renewal function of $H$.

By \cite[Theorem 5, page 79]{Ber} and \cite[Lemma 2.1]{KSV1}, $V$ is
absolutely continuous and has a continuous and strictly positive
density $v$ such that $v(0+)=1$.
The functions $V$ and $v$ enjoy the following estimates near the origin.

\begin{lemma}
{\rm(\cite[Lemma 2.2]{KSV1})} \label{l:estimate-for-V} Let $R>0$.
There exists a constant $c=c(R) >1$ such that for all $x\in (0,R]$,
we have $ c^{-1} \le v(x) \le c $ and $ c^{-1} x\le V(x) \le c x\, .
$
\end{lemma}

By \cite[Proposition 2.4]{KSV3}, the Laplace exponent $\chi$ of the
ladder height process $H$ of $Z_t$ is also  a complete Bernstein
function. Using this and the fact that $\chi$ has a drift,
 we see from \cite[Corollary 2.3]{KSV2}, that $v$ is completely
monotone. In particular, $v$ and the renewal function $V$ are
$C^{\infty}$ functions.

We will use $\bR^d_+$ to denote the half-space $\{ x=(x_1, \dots,
x_{d-1}, x_d):=(\tilde{x}, x_d)  \in \bR^d: x_d > 0 \}$.
 Define $w(x):=V((x_d)^+)$.

\begin{thm}\label{t:Sil}
The function $w$ is harmonic in $\R^d_+$ with respect to $X$ and,
for any $r>0$, regular harmonic in $\R^{d-1}\times (0, r)$ with
respect to $X$.
\end{thm}

\pf Since $Z_t:=W^d_{S_t}$ has a transition density, it satisfies
the condition ACC in  \cite{Sil}, namely the resolvent kernels are
absolutely continuous. The assumption in  \cite{Sil} that $0$ is
regular for $(0,\infty)$ is also satisfied since $X$ is of unbounded
variation. Further, by symmetry of $Z$, the notions of
coharmonic and harmonic functions coincide.
Now the theorem follows by the same argument as in \cite[Theorem 4.1]{KSV3}.
\qed

Unlike \cite[Proposition 4.2]{KSV3}, we prove the next result
without using the boundary Harnack principle.

\begin{prop}\label{c:cforI}
For all positive constants $r_0$ and $L$, we have
$$
\sup_{x \in \R^d:\, 0<x_d <L} \int_{B(x, r_0)^c \cap \bR^d_+}
w(y) j(|x-y|)\, dy < \infty\, .
$$
\end{prop}

\pf Without loss of generality, we assume $\wt x=0$.
We consider two separate cases.

\noindent (a) Suppose $L>x_d \ge r_0/4$.  By
\eqref{e:levy} and Theorem \ref{t:Sil},
for every $x \in \R^d_+$,
\begin{eqnarray}
w(x)
&\ge& \E_x\left[w\big(X_{\tau_{ B(x, r_0/2)\cap \bR^d_+}}\big):
X_{\tau_{ B(x, r_0/2)\cap \bR^d_+}} \in  {B(x, r_0)}^c \cap
\bR^d_+  \right]\nn\\
&=& \E_x \left[\int_0^{\tau_{B(x, r_0/2)\cap \bR^d_+}}   \int_{{B(x,
r_0)}^c \cap \bR^d_+} j(|X_t-y|)w(y)\, dydt \right]\, .
\label{e:fdfsaff}
\end{eqnarray}
Since $|z-y| \le |x-z|+|x-y| \le r_0 +|x-y| \le 2|x-y|$ for $(z,y)
\in B(x, r_0/2) \times B(x, r_0)^c$, using \eqref{H:1} and
\eqref{H:2}, we have $j(|z-y|) \ge c_1 j(|x-y|)$.
Thus, combining this with \eqref{e:fdfsaff}, we obtain that
$$
\int_{B(x, r_0)^c \cap \bR^d_+} w(y) j(|x-y|) dy  \le
c_1^{-1} \frac{w(x)}{\E_{x}[\tau_{B(x, r_0/2)\cap \bR^d_+}]} \le c_1^{-1} \frac{V(L)}{\E_{0}[\tau_{B(0, r_0/4)}]}.
$$

\noindent (b) Suppose $x_d  < r_0/4$. Note that if $|y-x|>r_0$, then
$|y|\ge |y-x|-|x| >3r_0/4$ and $|y| \le |y-x| +|x| \le |y-x|+ r_0/4
\le |y-x|+ |y-x|/4$. Thus, using \eqref{H:1} and \eqref{H:2}, we
have $j(|y|) \ge c_2 j(|x-y|)$ and
\begin{eqnarray}
\sup_{x \in \R^d:\, 0<x_d < r_0/4} \int_{B(x, r_0)^c \cap \bR^d_+}
w(y) j(|x-y|) dy  \le c_3 \int_{B(0, r_0/2)^c \cap \bR^d_+} w(y)
j(|y|) dy . \label{e:fdfsaff0}
\end{eqnarray}
Let $x_1:=(\wt 0, r_0/8)$.
By Theorem \ref{t:Sil} and \eqref{e:levy},
\begin{eqnarray}
\infty >w(x_1)
&\ge& \E_{x_1}\left[w(X_{\tau_{ B(0, r_0/4)\cap \bR^d_+}}):
X_{\tau_{ B(0, r_0/4)\cap \bR^d_+}} \in  {B(x, r_0/2)}^c \cap
\bR^d_+  \right]\nn\\
&=&
\E_{x_1} \left[\int_0^{\tau_{B(0, r_0/4)\cap \bR^d_+}}
 \int_{{B(0, r_0/2)}^c \cap \bR^d_+} j(|X_t-y|)w(y)\, dy\,dt \right] .
\label{e:fdfsaff2}
\end{eqnarray}
Since $|z-y| \le |z|+|y| \le (r_0/4) +|y| \le 2|y|$ for $(z,y)
\in B(0, r_0/4) \times B(0, r_0/2)^c$, using \eqref{H:1} and
\eqref{H:2}, we have $j(|z-y|) \ge c_3 j(|y|)$.
Thus, combining this with \eqref{e:fdfsaff2}, we obtain that
\begin{eqnarray}
\infty >w(x_1) &>& c_3  \E_{x_1} \left[\int_0^{\tau_{B(0, r_0/4)
\cap \bR^d_+}}   \int_{{B(0, r_0/2)}^c \cap \bR^d_+}  j(|y|)w(y)\,
dy\,dt \right]\nn\\
&=&
c_3 \E_{x_1}[\tau_{B(0, r_0/4)\cap \bR^d_+}] \int_{B(0, r_0/2)^c
\cap \bR^d_+}  j(|y|)w(y) dy.
\label{e:fdfsaff3}
\end{eqnarray}
Combining \eqref{e:fdfsaff0} and \eqref{e:fdfsaff3}, we conclude
that
$$
\sup_{x \in \R^d:\, 0<x_d < r_0/4}
\int_{B(x, r_0)^c \cap \bR^d_+} w(y) j(|x-y|) dy \le
c_4\frac{V(r_0/8)}{\E_{0}[\tau_{B(0, r_0/8)}] } < \infty.
$$
\qed

We now define an operator ($\Delta+ \sA$, $\mathfrak{D}
(\Delta+\sA)$) as follows:
\begin{align}
\sA  f(x)
:=&\lim_{\eps \downarrow 0}
\int_{B(x, \eps)^c}
\left(f(y)-f(x)\right)j(|y-x|)\, dy, \nn\\
\mathfrak{D}(\Delta+\sA):=&\left\{f \in C^2(\R^d): \lim_{\eps \downarrow 0}
\int_{B(x, \eps)^c}
\left(f(y)-f(x)\right)j(|y-x|)\, dy \text{ exists and is finite } \right\}.
\label{generator}
\end{align}
Recall that $C^2_0(\bR^d)$ is the collection of $C^2$ functions in $\RR^d$ vanishing at
infinity
It is well known that
$C^2_0 (\bR^d)\subset \mathfrak{D}(\Delta+\sA)$ and that, by the rotational symmetry of $X$, $\Delta+\sA$
restricted to $C^2_0(\bR^d)$ coincides with the infinitesimal generator of
the process $X$ (e.g. \cite[Theorem 31.5]{Sa}).

The proof of the next result is similar to that of \cite[Theorem
4.3]{KSV3}. We give the proof here for completeness.
\begin{thm}\label{c:Aw=0}
$(\Delta+\sA) w(x)$ is well defined and $(\Delta+\sA) w(x)=0$ for
all $x \in \bR^d_+$.
\end{thm}

\pf It follows from Proposition \ref{c:cforI} and the fact that $j$
is a L\'evy density that for
any $L>0$ and $\varepsilon \in (0,1/2)$
\begin{eqnarray}
\lefteqn{\sup_{x \in \R^d:\ 0<x_d <L}
\left|\int_{B(x, \varepsilon)^c}
(w(y)-w(x)) j(|y-x|)dy \right|\nn}\\
&\le&
\sup_{x \in \R^d:\ 0<x_d <L}
\int_{B(x, \varepsilon)^c}
w(y)j(|y-x|)\, dy
+V(L)\int_{B(x, \varepsilon)^c}
j(|y|)dy
<\infty\, . \label{e:dedwer}
\end{eqnarray}
Hence, for every $\varepsilon \in (0,1/2)$,
$
\sA_\eps w(x):=\int_{B(x, \varepsilon)^c}(w(y)-w(x))j(|y-x|)dy
$
is well defined.
Using the smoothness of $w$ in $\R^d_+$ and following the same argument in
\cite[Theorem 4.3]{KSV3},  we can show that
$\sA w$ is well defined in $\R^d_+$ and
$\sA_\eps w(x)$ converges to
$$
\sA w(x)=\int_{\RR^d}\left(w(y)-w(x)-{\bf 1}_ {\{|y-x|<1\}} (y-x)
\cdot\nabla w(x)\right)j(|y-x|) dy
$$
locally uniformly in ${\R^d_+}$ as $\varepsilon\to 0$ and
the
function $\sA w(x)$ is continuous in ${\R^d_+}$.

Suppose that $U_1$ and $U_2$ are relatively compact open subsets of
$\RR^d_+$ such that $\overline{U_1} \subset U_2 \subset
\overline{U_2} \subset \RR^d_+$.
It follows again from the same argument in \cite[Theorem 4.3]{KSV3}
 that the conditions \cite[(2.4),
(2.6)]{C} are true.
 Thus, by \cite[Lemma 2.3, Theorem 2.11(ii)]{C},
we have that for any $f\in C^2_c(\R^d_+)$, \bee\label{e:C2.11}
0=\int_{\R^d} \nabla w(x) \cdot \nabla f(x) \, dx +
\frac12\int_{\R^d}\int_{\R^d}(w(y)-w(x))(f(y)-f(x))j(|y-x|)\, dx\,
dy. \eee For $f\in C^2_c(\R^d_+)$ with supp$(f) \subset
\overline{U_1} \subset U_2 \subset \overline{U_2} \subset \RR^d_+$,
\begin{align*}
&\int_{\R^d}\int_{\R^d}|w(y)-w(x)||f(y)-f(x)|j(|y-x|) dxdy\\
=&\int_{U_2}\int_{U_2}|w(y)-w(x)||f(y)-f(x)|j(|y-x|) dxdy+
2\int_{U_1}\int_{U_2^c}|w(y)-w(x)||f(x)|j(|y-x|) dxdy\\
\le&c_1\int_{U_2 \times U_2}
|y-x|^2 j(|y-x|) dxdy+2\|f\|_\infty|U_1|  \left(\sup_{x \in U_1}
w(x)\right) \int_{U_2^c}j(|y-x|) dy\\
&+2\|f\|_\infty\int_{U_1} \int_{U_2^c} w(y) j(|x-y|)dydx
\end{align*}
is finite by
Proposition \ref{c:cforI} and the fact that $j(|x|)dx$ is a L\'evy measure.
Thus by \eqref{e:C2.11}, Fubini's theorem and the dominated
convergence theorem, we have for any $f\in C^2_c(\R^d_+)$,
\begin{align*}
&0=\int_{\R^d} \nabla w(x) \cdot \nabla f(x) \, dx + \frac12
\lim_{\varepsilon\downarrow0}\int_{\{(x, y)\in {\R^d}\times {\R^d},\
|y-x|>\varepsilon\}}(w(y)-w(x))(f(y)-f(x))j(|y-x|)\, dx\, dy\\
&=-\int_{\R^d} \Delta w(x)  f(x) \, dx
-\lim_{\varepsilon\downarrow0}\int_{\R^d_+} f(x)
\left(\int_{B(x, \varepsilon)^c}(w(y)-w(x))
j(|y-x|) dy \right) dx\\
&=- \int_{\R^d} \Delta w(x)  f(x) \, dx  -\,\int_{\R^d_+} f(x)\sA
w(x)\, dx\, =\,- \int_{\R^d}( \Delta+\sA) w(x)  f(x) \, dx
\end{align*}
where we have used the fact $\sA_\eps w \to \sA w$ converges
uniformly on the support of $f$. Hence, by the continuity of $(
\Delta+\sA) w$, we have $( \Delta+\sA) w(x)=0$ in ${\R^d_+}$. \qed

\section{Analysis on $C^{1,1}$ open set}

Recall that  $\Lambda \ge 1$ and that $D$ is a $C^{1,1}$ open set
with characteristics $(R, \Lambda)$ and $D$ satisfies the  uniform
interior ball condition and the uniform exterior ball condition with
radius $R \le 1$.
The proof of the next lemma is motivated by that of \cite[Lemma 4.4]{KSV3}.

\begin{lemma}\label{L:Main}
Fix $Q \in \partial D$ and define
$$
h(y):=V\left(\delta_D (y)\right){\bf 1}_{D\cap B(Q, R)}(y).
$$
There exists $C_1=C_1(
\Lambda, R)>0$  independent of
 $Q$ such that $( \Delta+\sA) h$ is well
defined in $D\cap B(Q, R/4)$ a.e. and
\bee\label{e:h3}
|(\Delta+\sA) h(x)|\le C_1 \quad \text{for  a.e. } x \in D\cap B(Q, R/4)\, .
\eee
\end{lemma}

\pf In this proof, we fix $x \in D\cap B(Q, R/4)$ and
$x_0\in\partial D$ satisfying $\delta_D(x)=|x-x_0|$.  We also fix
the $C^{1, 1}$ function $\varphi$ and the coordinate system
$CS=CS_{x_0}$ in the definition of $C^{1, 1}$ open set so that $x=(
0, x_d)$ with $0<x_d <R/4$ and $B(x_0, R)\cap D=\{ y=(\wt y, \, y_d)
\in B(0, R) \mbox{ in } CS : y_d > \varphi (\wt y) \}.$ Let
$$
\varphi_1(\wt y):=R -\sqrt{ R^2-|\wt y|^2}
\quad \text{and}\quad \varphi_2(\wt y):=-R +\sqrt{ R^2-|\wt
y|^2}.
$$
Due to the  uniform interior ball condition and the uniform
exterior ball condition with radius $R$, we have
 \bee\label{e:phi012}
\varphi_2(\wt y)  \le \varphi (\wt y) \le \varphi_1(\wt y) \quad \text{for
every } y \in D\cap B(x, R/4).
 \eee
Define $H^+:= \left\{y=(\wt y, \, y_d) \in CS:y_d>0 \right\}$ and
let
$$
A:=\{y=(\widetilde{y},y_d) \in (D \cup H^+)\cap B(x, R/4):
\varphi_2(\wt y)  \le y_d  \le \varphi_1(\wt y)\},
$$
$$
E:=\{y=(\widetilde{y},y_d) \in B(x, R/4):  y_d  > \varphi_1(\wt
y)\}.
$$
Note that, since $|y-Q| \le |y-x| +|x-Q| \le R/2$ for $y
\in B(x, R/4)$, we have $B(x,  R/4) \cap D \subset
B(Q,   R/2)\cap D.$

Let
$$
h_{x}(y):=V\left(\delta_{_{H^+}}(y)\right).
$$
Note that $h_{x}(x)=h(x)$. Moreover, since
$\delta_{_{H^+}}(y)=(y_d)^+$ in $CS$,
it follows from Theorem \ref{c:Aw=0} that $\sA h_x$ is well defined in $H^+$ and
\bee\label{e:hz}
(\Delta+\sA)  h_{x}(y)=0, \quad  \forall y\in H^+.
\eee
We show now that $\sA (h-h_x)(x)$ is well defined.
For each $\varepsilon >0$ we have that
\begin{align*}
&\bigg|\int_{\{y \in D \cup H^+: \, |y-x|>\varepsilon\}}
{(h(y)-h_{x}(y))}j(|y-x|)\ dy\bigg|
\nn\\
\leq & \int_{B(x,R/4)^c}
(h(y)+h_{x}(y))
j(|y-x|)dy +\int_{A} (h(y)+h_{x}(y))j(|y-x|)\ dy
\\&\quad +\int_{E}
{|h(y)-h_{x}(y)|}j(|y-x|) dy
\,=:\,I_1\,+\,I_2\,+\,I_3.
\end{align*}

By the fact that $h(y)=0$ for $y \in B(Q, R)^c$,
\begin{align}
I_1 \le &\sup_{z \in \R^d:\ 0<z_d <R} \int_{B(z,  R/4)^c \cap H^+}
V(y_d) j(|z-y|) dy+c_1\int_{B(0,  R/4)^c}j (|y|)dy =:K_1+K_2.
\nonumber
\end{align}
$K_2$ is clearly finite since $J$ is the L\'evy density of $X$ while
$K_1$ is  finite by Proposition \ref{c:cforI}.

For $y \in A$, since $V$ is increasing and $(R - \sqrt{R^2-|\wt
y|^2}) \le  R^{-1}|\wt y|^2$, we see that
\begin{align}
h_{x}(y)+h(y) \le 2V(\varphi_1 ( \wt y) -\varphi_2 ( \wt y)) \le2
V(2 R^{-1}|\wt y|^2) \le 2V(2 R^{-1}|y-x|^2). \label{e:dfe2}
\end{align}
Using \eqref{e:dfe2} and Lemma  \ref{l:estimate-for-V}, we have
\begin{align}
I_2\leq&c_2 \int_{ A}   |y-x|^2  j(|y-x| )dy  \le c_2  \int_{ B(0,
R/4)}   |z|^2  j(|z| )dz < \infty . \label{e:fgy6} \end{align}

For $I_3$, we consider two cases separately: If $0
<y_d=\delta_{_{H^+}}({y}) \le  \delta_D({y})$, since $v$ is
decreasing,
\begin{align}
h(y)-h_{x}(y) \le V(y_d+R^{-1}|\wt y|^2) -V(y_d) =
\int_{y_d}^{y_d+R^{-1}|\wt y|^2} v(z)dz \le R^{-1}|\wt y|^2
v(y_d).  \label{e:KG}
\end{align}
If $y_d=\delta_{_{H^+}}({y}) >  \delta_D({y})$ and $y \in E$, using
the fact that $\delta_D({y})$ is greater than or equal to the
distance between $y$ and the graph of $\varphi_1$ and
\begin{align*}
y_d-R+\sqrt{ |\wt y|^2+(R-y_d)^2} &= \frac{|\wt y|^2} {\sqrt{ |\wt
y|^2+(R-y_d)^2} + (R-y_d)}\,\le\, \frac{ |y-x|^2} {2 (R-y_d)} \le
\frac{ |y-x|^2} {R},
\end{align*}
we have
\begin{align}
h_{x}(y)-h(y) \le\int^{y_d}_{R-\sqrt{ |\wt y|^2+(R-y_d)^2}}
v(z)dz\le  R^{-1}  |y-x|^2 \,v(R-\sqrt{ |\wt y|^2+(R-y_d)^2}).
\label{e:KG2}
\end{align}
Thus, by \eqref{e:KG}-\eqref{e:KG2} and Lemma \ref{l:estimate-for-V},
\begin{align*}
I_3 \,\le\, c_3\int_{E} | y-x|^2 j(|y-x|) dy \le c_3  \int_{ B(0,
R/4)}   |z|^2  j(|z| )dz < \infty .
\end{align*}
We have proved
 \bee\label{e:I1234}
|\sA(h-h_x)(x)|\le I_1+I_2+I_3 \leq c_4
 \eee
for some constant $c_4=c_4(R, \Lambda)>0$.

The estimate \eqref{e:I1234} shows in particular that the limit
$$
\lim_{\varepsilon\downarrow 0}\int_{\{y \in D \cup H^+:
|y-x|>\varepsilon\}}{(h(y)-h_{x}(y))}j(|y-x|)\, dy
$$
exists and hence $\sA(h-h_x)(x)$ is well defined.

We now consider $\Delta(h-h_x)$.
Note that for  a.e.~$x \in D\cap B(Q, R/4)$, the second order
partial derivatives of the function $y \to \delta_D (y)$ exist at
$x$. Without loss of generality we assume that $x$ has been chosen
so that the second order partial derivatives of the function $y \to
\delta_D (y)$ exist at $x$.
 Since $h_{x}(y)=V\left((y_d)^+\right)$ in $CS$, we have $\Delta h_x(x)= v'(x_d)$.
 Moreover, since $\delta_D(y)=y_d$ for $y=(\wt 0, x_d+\eps)$
when $|\eps|$ is small, $\partial^2_{x_d} h(x)= v'(x_d)$.
Thus
\begin{align}\label{e:dsnew1}
&\Delta(h-h_x)(x) = \sum_{i=1}^{d-1}
\frac{\partial^2  V\left(\delta_D (y)\right)}{\partial y_i^2}|_{y=x}
=\sum_{i=1}^{d-1}
\frac{\partial}{\partial y_i} \left(v(\delta_D (y))
\frac{\partial \delta_D (y)}{\partial y_i}\right)|_{y=x} \nonumber \\
&= \sum_{i=1}^{d-1} v'(\delta_D (x))\left(\frac{\partial \delta_D (y)}{\partial y_i}|_{y=x} \right)^2
+ v(\delta_D (x))\frac{\partial^2 \delta_D (y)}{\partial y_i^2}|_{y=x}\ .
\end{align}
In the coordinate system $CS$,
\begin{equation}\label{e:dsnew2}
\frac{\partial \delta_D(y)}{\partial y_i}\large|_{y=x}=0 \quad \text{ and } \quad
\left|\frac{\partial^2\delta_D(y)}{\partial y_i^2}\large|_{y=x}
\right|\le \frac{4}{3R}, \quad i=1, \dots, d-1.
\end{equation}
Indeed, let $\epsilon \in \bR$ with $|\epsilon|$ small,
and $x_{\epsilon,i}:=(0, \dots, \epsilon, \dots 0, x_d)$, $i=1, \dots , d-1$.
Due to the  uniform interior ball condition and the uniform
exterior ball condition with radius $R$, we have
$$
R-\sqrt{\epsilon^2+(R-x_d)^2}-x_d \le \delta_D(x_{\epsilon,i})-\delta_D(x) \le \sqrt{\epsilon^2+(R+x_d)^2}-R-x_d\, ,
$$
so
\begin{align*}
\frac{1}{\epsilon}| \delta_D(x_{\epsilon,i})-\delta_D(x)|
\le& \frac{1}{\epsilon} \left( \sqrt{\epsilon^2+(R-x_d)^2}-(R-x_d)   \right) \vee  \frac{1}{\epsilon}\left(  \sqrt{\epsilon^2+(R+x_d)^2}-(R+x_d) \right)\\
=&\left(\frac{\epsilon}{\sqrt{\epsilon^2+(R-x_d)^2}+(R-x_d) }  \right) \vee  \left( \frac{\epsilon}{\sqrt{\epsilon^2+(R+x_d)^2}+(R+x_d)}  \right),
\end{align*}
which goes to zero as $\epsilon \to 0$. The bound involving the second partial derivatives can be proved in
a similar way using the elementary fact that $\frac{\partial^2\delta_D(y)}{\partial x_i^2}|_{y=x}
= \lim_{\epsilon \to 0} \frac{1}{h^2}( \delta_D(x_{\epsilon,i})+\delta_D(x_{-\epsilon,i})-2\delta_D(x)).$
Therefore, combining \eqref{e:dsnew1}, \eqref{e:dsnew2} and Lemma \ref{l:estimate-for-V}, we have
$$
|\Delta(h-h_x)(x)| \le  c_5\sum_{i=1}^{d-1}
\left|\frac{\partial^2\delta_D(y)}{\partial x_i^2}\large|_{y=x}
\right| \le c_6.
$$
Using this, \eqref{e:hz}, \eqref{e:I1234}, and linearity we get
that $(\Delta+\sA) h(x)$ is well defined and $|(\Delta+\sA)
h(x)|\le c_7$. \qed

We use $C^{\infty}_c(\R^d)$ to denote the space of infinitely
differentiable functions with compact support.
Using the fact that $\Delta+\sA$ restricted to
$C^{\infty}_c(\R^d)$  coincides with the
infinitesimal generator of the process $X$, we see that the
following Dynkin's formula is true: for $f \in C_c^{\infty}(\R^d)$ and
any bounded open subset $U$ of $\R^d$,
\begin{equation} \label{e:*334}
\E_x\int_0^{\tau_U}  {(\Delta+\sA)} f(X_t) dt
=\E_x[f(X_{\tau_U})]- f(x).
\end{equation}

\begin{lemma}\label{l2.1}
For every $r_1> 0$ and every  $a \in (0,1)$, there exists a
positive constant $c=c(r_1, d, a)$ such that for any $r\in
(0, r_1]$ and any open sets $U$ and $D$ with $B(0, a r ) \cap D
\subset U  \subset  D$, we have
$$
\P_x\left(X_{\tau_U} \in D\right) \,\le\, c\,r^{-2}\, \E_x[\tau_U],
\qquad x \in D\cap B(0, a r/2).
$$
\end{lemma}

\pf For fixed $a \in (0,1)$, take a sequence of radial functions
$\phi_m$ in $C^{\infty}_c(\R^d)$ such that $0\le \phi_m\le 1$,
$$
\phi_m(y)=\left\{
\begin{array}{lll}
0, & |y|< a/2\\
1, & a\le |y|\le m+1\\
0, & |y|>m+2,
\end{array}
\right.
$$
and that $\|\sum_{i, j}|\frac{\partial^2}{\partial y_i\partial y_j}
\phi_m|\|_\infty < c_1=c_1(a)< \infty$.
Define $\phi_{m, r}(y)=\phi_m(\frac{y}{r})$ so that
$0\le \phi_{m, r}\le 1$,
\begin{equation}\label{e:2.11}\phi_{m, r}(y)=
\begin{cases}
0, & |y|<ar/2\\
1, & ar\le |y|\le r(m+1)\\
0, & |y|>r(m+2),
\end{cases}
\quad \text{and} \quad  \sup_{y\in \R^d} \sum_{i,
j}\left|\frac{\partial^2}{\partial y_i\partial y_j} \phi_{m,
r}(y)\right| \,<\, c_1\, r^{-2}.
\end{equation}
Using \eqref{e:2.11}, we see that
\begin{eqnarray}
&&  \left|\int_{\R^d} (\phi_{m,r}(x+y)-\phi_{m,r}(x)-(\nabla
\phi_{m,r}(x) \cdot y)1_{B(0, 1)}(y))J(y) dy \right|\nn\\
&&\le \left|\int_{\{|y|\le 1\}}
(\phi_{m,r}(x+y)-\phi_{m,r}(x)-(\nabla \phi_{m,r}(x) \cdot y)1_{B(0, 1)}(y))J(y) dy\right|+
2\int_{\{|y|>1\}}J(y) dy \nn\\
&&\le \frac{c_2}{r^2}\int_{\{|y|\le 1 \}} |y|^2 J(y)dy + 2
\int_{\{|y|>1\}}J(y) dy  \le \frac{c_3}{r^2}\, . \label{e2.100}
\end{eqnarray}
Now, by combining \eqref{e:*334}, (\ref{e:2.11}) and (\ref{e2.100}),
we get that for any $x\in D\cap B(0, ar/2)$,
\begin{align*}
&\P_x\left(X_{\tau_U} \in \{ y \in D: ar \le |y| <(m+1)r \}\right)
=\E_x\left[\phi_{m, r} \left(X_{\tau_U}\right): X_{\tau_U} \in \{ y
\in D: ar  \le |y| <(m+1)r
\}\right]\\
&\le \E_x\left[\phi_{m, r} \left(X_{\tau_U}\right)\right]=
\E_x\left[ \int_0^{\tau_U}   {(\Delta+\sA)} \phi_{m, r}(X_t)dt \right] \le
\|{(\Delta+\sA)}\phi_{m, r}\|_\infty\,  \E_x[\tau_U]
\le c_{4} r^{-2}\E_x[\tau_U].
\end{align*}
Therefore, since $B(0, a r ) \cap D \subset U$,
\begin{align*}
\P_x\left(X_{\tau_U} \in D\right)&=\P_x\left(X_{\tau_U} \in
\{ y \in D: ar \le |y| \}\right) \\
&= \lim_{m\to
\infty}\P_x\left(X_{\tau_U} \in \{ y \in D: a r \le |y| <(m+1)r
\}\right) \,\le\,
c_5\,r^{-2}\E_x[\tau_U].
\end{align*}
\qed

Define $\rho_Q (x) := x_d -  \varphi_Q (\wt x),$ where $(\wt x,
x_d)$ are the coordinates of $x$ in $CS_Q$. Note that for every $Q
\in \partial D$ and $ x \in B(Q, R)\cap D$ we have
\begin{equation}\label{e:d_com}
(1+\Lambda^2)^{-1/2} \,\rho_Q (x) \,\le\, \delta_D(x)  \,\le\,
\rho_Q(x).
\end{equation}
We define for $r_1, r_2>0$
$$
D_Q( r_1, r_2) :=\left\{ y\in D: r_1 >\rho_Q(y) >0,\, |\wt y | < r_2
\right\}.
$$

Let $R_1 :=R/(4\sqrt{1+(1+ \Lambda)^2})$.
By Lemma \ref{l:estimate-for-V}, $V(\delta_D(x)$ on the right-hand sides of
\eqref{e:L:2}--\eqref{e:L:3} can be replaced by $\delta_D(x)$. The reason
we prefer the forms below is that the function $V$ will be used in the proof.

\begin{lemma}\label{L:2}
There are constants $\lambda_0 >2 R_1^{-1} $, $\kappa_0  \in (0,1)$ and $c=c(R,  \Lambda
)>0$ such that for every $\lambda \ge \lambda_0$, $Q \in
\partial D$ and $x \in D_Q (2^{-1}(1+\Lambda)^{-1}\kappa_0 \lambda^{-1} ,
\kappa_0\lambda^{-1} )$ with $\wt x=0$,
  \bee\label{e:L:2}
\P_{x}\left(X_{ \tau_{ D_Q (   \kappa_0 \lambda^{-1} , \lambda^{-1} )}} \in
D\right) \,\le\, c \, \lambda \, V( \delta_D (x))
 \eee
and
 \bee\label{e:L:3}
\E_x\left[ \tau_{ D_Q  (   \kappa_0 \lambda^{-1} , \lambda^{-1}  )}
\right]\,\le\, c\, \lambda^{-1}\, V(\delta_D (x)).
 \eee
\end{lemma}

\pf Without loss of generality, we assume $Q=0$. Let
$\varphi=\varphi_0:\bR^{d-1}\to \bR$ be the $C^{1,1}$ function  and
$CS_0$ be the coordinate system in the definition of $C^{1, 1}$ open
set so that $B(0, R)\cap D= \big\{(\wt y, \, y_d) \in B(0,
R)\textrm{ in } CS_0: y_d > \varphi (\wt y) \big\}.$ Let $\rho (y)
:= y_d -  \varphi (\wt y)$ and $D ( a, b):=D_0 ( a, b)$.

Note that
 \bee\label{e:lsd}
|y|^2 = |\wt y|^2+ |y_d|^2 <r^2 +(|y_d- \varphi(\wt y)|+ |\varphi(\wt
y)|)^2 < (1+(1+ \Lambda)^2) r^2 \quad \text{for every } y \in
D(r,r)\, .
\eee
By this and the definition of $R_1$, we
have $D ( r, s) \subset D(R_1, R_1)\subset B(0, R/4) \cap D \subset
B(0, R)\cap D$  for every $ r,s \le R_1 . $

Using Lemma \ref{l:estimate-for-V}, we can and will
choose $\delta_0 \in (0, R_1)$ small such that
 $$
2 r^2 \le V( (1+ \Lambda^2)^{-1/2} r) \quad \text{ for all }
r \le 4\delta_0.
$$
Then, by \eqref{e:d_com}, the subadditivity and monotonicity of $V$,
for every $\lambda \ge 1$ and every $ y \in B(0,R) \cap D$ with
$\rho (y) \le 4\lambda^{-1}\delta_0$, we have
\begin{equation}\label{e:ssp1}
2 \lambda^2 \rho (y)^2 \le V(\lambda \delta_D(y))  \le (\lambda+1)
V( \delta_D(y)) \le 2\lambda V( \delta_D(y)).
\end{equation}

Since $\Delta \varphi(\tilde{y})$ is well defined for a.e.~$\tilde{y}$
with respect to the $(d-1)$-dimensional Lebesgue measure, it follows that
$\Delta \rho(y)$ exists for a.e.~$y$ with respect to the
$d$-dimensional Lebesgue measure.
Using the fact that the derivative
of a Lipschitz function is essentially bounded by its Lipschitz
constant, we have
for a.e.~$y\in B(0,R)\cap D$ that
\begin{eqnarray*}
\Delta \rho (y)^2 = \Delta (y_d -  \varphi (\wt y))^2 =2(1+ |\nabla
\varphi (\wt y)|^2) - 2\, \rho (y) \Delta \varphi(\wt y) \ge  2(1-
\rho (y) \| \Delta \varphi\|_\infty)\,  .
\end{eqnarray*}
Choosing $\delta_0 \in (0, R_1)$ smaller if necessary we can get that
\begin{equation}\label{e:ssp2}
\Delta \rho (y)^2 \ge 1  \quad \text{ for a.e. } y \in B(0,R) \cap D
\text{ with }\rho (y) \le 2\delta_0.
\end{equation}

Let $g(y)=g(\wt y, y_d)$ be a smooth  function on $\bR^d$
with $0 \le g(\wt y, y_d) \le 2$, $g(\wt y, y_d) \le y_d^2$,
\begin{equation}\label{e:ssp101}
\sum_{i, j=1}^d|\frac{\partial^2 g}{\partial y_i\partial y_j} |+\sum_{i=1}^d
|\frac{\partial g}{\partial y_i} | \le c_1,
\end{equation}
and
$$
g(y)=
\begin{cases}
0, & \text{ if }  -\infty < y_d <0, \text{ or } y_d \ge 4 \text{ or } |\wt y|>2 \\
y_d^2, & \text{ if }  0\le y_d<1 \text{ and } |\wt y|<1   \\
-(y_d-2)^2+2, & \text{ if }  1 \le y_d \le 3 \text{ and } |\wt y|<1\\
(y_d-4)^2, & \text{ if }  3 \le y_d \le 4 \text{ and } |\wt y|<1.
\end{cases}
$$
Thus $\mathrm{supp}(g) \subset \{   |\wt y| \le 2, 0 \le y_d \le 4\}$.

For $\lambda  >1$, let
$
g_\lambda(y):=g_\lambda(\wt y, y_d):=g( \lambda \delta_0^{-1}
 \wt y,  \lambda \delta_0^{-1} \rho(y))
$
so that
\begin{equation}\label{e:new001}
\mathrm{supp}
(g_\lambda) \subset \{   |\wt y| \le 2\lambda^{-1} \delta_0, \ \, 0 \le \rho(y)  \le 4\lambda^{-1} \delta_0 \}.
\end{equation}
Then, since $\sum_{i, j=1}^d|\frac{\partial^2}{\partial y_i\partial y_j} g(y)|$
is essentially bounded, using \eqref{e:ssp101}, we have
\begin{equation}\label{e:ssp3}
\sum_{i, j=1}^d|\frac{\partial^2}{\partial y_i\partial y_j}\, g_\lambda
(y)| \le c_2 \lambda^2 \quad \text{a.e. } y.
\end{equation}
Note that, by the definition of $g$, $g_\lambda(y)=\lambda^2
 \delta_0^{-2} \rho(y)^2$ on $D(\lambda^{-1} \delta_0, \lambda^{-1}
\delta_0)$. Thus,  from \eqref{e:ssp2} we get
\begin{equation}\label{e:ssp4}
\Delta g_\lambda(y) \ge \lambda^2 \delta_0^{-2}  \quad \text{ for
a.e. } y \in D(\lambda^{-1} \delta_0, \lambda^{-1} \delta_0).
\end{equation}

On the other hand, by \eqref{e:ssp3} we have
\begin{eqnarray*}
&&  \left|\int_{\R^d} (g_{\lambda}(y+z)-g_{\lambda}(y)-(\nabla
g_{\lambda}(y) \cdot z)1_{B(0, \lambda^{-1} )}(z))J(z)\, dz \right| \nonumber\\
&&\le \left|\int_{\{|z|\le \lambda^{-1} \}}
(g_{\lambda}(y+z)-g_{\lambda}(y)-(\nabla g_{\lambda}(y) \cdot
z)1_{B(0, \lambda^{-1} )}(z))J(z)\, dz \right| \nonumber \\
&&\quad+\int_{\{\lambda^{-1} <|z| \le 1\}}J(z) g_{\lambda}(y+z)dz+
\left(\int_{\{\lambda^{-1} <|z| \le 1\}}J(z) dz\right)g_{\lambda}(y)
  +2\int_{\{1 <|z|\}}J(z)\, dz
\nonumber\\
&&\le c_3\lambda^2\int_{\{|z|\le \lambda^{-1}  \}}
|z|^2 J(z)\, dz+
2\int_{\{1 <|z|\}}J(z)\, dz
\nonumber
\\
&&\quad+
\int_{\{\lambda^{-1} <|z| \le 1\}}J(z) g_{\lambda}(y+z)dz
+\left(\int_{\{\lambda^{-1} <|z| \le 1\}}J(z) dz\right)g_{\lambda}(y).
\end{eqnarray*}
Thus
\begin{eqnarray}\label{e:ssp51}
&&\lambda^{-2}
\left|\int_{\R^d} (g_{\lambda}(y+z)-g_{\lambda}(y)-(\nabla
g_{\lambda}(y) \cdot z)1_{B(0,\lambda^{-1})}(z))J(z)
 dz \right| \nonumber \\
 &\le& c_3\int_{\{|z|\le \lambda^{-1}  \}}
|z|^2 J(z)\, dz+
2\lambda^{-2}\int_{\{1 <|z|\}}J(z)\, dz
\nonumber
\\
&&\quad+\lambda^{-2}
\int_{\{\lambda^{-1} <|z| \le 1\}}J(z) g_{\lambda}(y+z)dz
+\lambda^{-2}\left(\int_{\{\lambda^{-1} <|z| \le 1\}}J(z) dz\right)g_{\lambda}(y)\nonumber\\
 &\le& c_3\int_{\{|z|\le \lambda^{-1}  \}}
|z|^2 J(z)\, dz+
2\lambda^{-2}\int_{\{1 <|z|\}}J(z)\, dz
\nonumber
\\
&&\quad+
\int_{\{\lambda^{-1} <|z| \le 1\}}J(z)|z|^2 g_{\lambda}(y+z)dz
+\left(\int_{\{0 <|z| \le 1\}}|z|^2J(z) dz\right)g_{\lambda}(y).
\end{eqnarray}

We claim that for every $\lambda  >1$ and $y \in D(\lambda^{-1} \delta_0, \lambda^{-1} \delta_0)$, the function $z \to g_{\lambda}(y+z)$ is supported in
$B(0,
3 \lambda^{-1}\delta_0 \sqrt{(4\Lambda)^2+1})$.

Fix $\lambda  >1$ and $y \in D(\lambda^{-1} \delta_0, \lambda^{-1} \delta_0)$ and
suppose that $z \in  B(0, 3 \lambda^{-1}\delta_0 \sqrt{(4\Lambda)^2+1})^c$. Then
either $|\wt z| \ge 3 \lambda^{-1}\delta_0$, or   $|\wt z| < 3 \lambda^{-1}\delta_0$ and $| z_d| \ge 12 \lambda^{-1}\delta_0 \Lambda$.
If $|\wt z| \ge 3 \lambda^{-1}\delta_0$, then clearly
$|\wt y + \wt z| \ge |\wt z|-|\wt y| \ge 3 \lambda^{-1}\delta_0-  \lambda^{-1}\delta_0 = 2 \lambda^{-1}\delta_0.$ Thus by  \eqref{e:new001},
$g_{\lambda}(y+z)=0$.
Now assume $|\wt z| < 3 \lambda^{-1}\delta_0$ and $|z_d| \ge 12 \lambda^{-1}\delta_0 \Lambda$.
If $z_d \le -12 \lambda^{-1}\delta_0 \Lambda$, then $g_{\lambda}(y+z)=0$.
If $z_d \ge 12 \lambda^{-1}\delta_0 \Lambda$, we have
$$
\rho(y+z) \ge z_d -|\psi(\wt y+ \wt z)|  \ge 12 \lambda^{-1}\delta_0 \Lambda - \Lambda ( |\wt z|+|\wt y|)
 \ge \lambda^{-1} \Lambda (12 \delta_0 - 3\delta_0- \delta_0) = 7 \lambda^{-1} \Lambda \delta_0.
$$ Thus by  \eqref{e:new001},
$g_{\lambda}(y+z)=0$.
The claim is proved.

Using the above claim and the fact that $g_\lambda(y)=\lambda^2
 \delta_0^{-2} \rho(y)^2$ on $D(\lambda^{-1} \delta_0, \lambda^{-1}
\delta_0)$,
we have from \eqref{e:ssp51}, that for $y \in D(\lambda^{-1} \delta_0, \lambda^{-1} \delta_0)$
\begin{align}\label{e:ssp52}
&\lambda^{-2}
\left|\int_{\R^d} (g_{\lambda}(y+z)-g_{\lambda}(y)-(\nabla
g_{\lambda}(y) \cdot z)1_{B(0,\lambda^{-1})}(z))J(z)
 dz \right| \nonumber \\
 \le& c_3\int_{\{|z|\le \lambda^{-1}  \}}
|z|^2 J(z)\, dz+
2\lambda^{-2}\int_{\{1 <|z|\}}J(z)\, dz
\nonumber
\\
&\quad+
\int_{\{\lambda^{-1} <|z| \le 1 \wedge
3 \lambda^{-1}\delta_0 \sqrt{(4\Lambda)^2+1}\}}J(z)|z|^2 dz
+c_4\lambda^2
 \delta_0^{-2} \rho(y)^2\nonumber \\
  \le& (c_3+1) \int_{\{|z|\le
  3 \lambda^{-1}\delta_0 \sqrt{(4\Lambda)^2+1} \}}
(1 \wedge |z|^2) J(z)\, dz+
2\lambda^{-2}\int_{\{1 <|z|\}}J(z)\, dz
+c_4\lambda^2
 \delta_0^{-2} \rho(y)^2,
\end{align}
where $c_4:=
2^{-1} \vee \int_{\{0 <|z| \le 1\}}|z|^2J(z) dz$.
 Define
$$
h(y) := V(\delta_D(y)) {\bf 1}_{B(0, R) \cap D}(y) \quad \text{ and}
\quad h_\lambda(y) := \lambda h(y) - g_{\lambda} (y).
$$
Choose $\lambda_* \ge 2
$ large such that  for every $\lambda \ge
\lambda_*$,
$$
 (c_3+1) \int_{\{|z|\le 2 \lambda^{-1}\delta_0 \sqrt{(4\Lambda)^2+1} \}}
(1 \wedge |z|^2) J(z)\, dz+
2\lambda^{-2}\int_{\{1 <|z|\}}J(z)\, dz  \le  4^{-1} \delta_0^{-2}\,\quad
\text{ and } \quad  \frac{1}{4} \lambda \delta_0^{-2} \ge  C_1,
$$
where $C_1$ is the constant from Lemma \ref{L:Main}.
Then by
\eqref{e:ssp4} and \eqref{e:ssp52},   for every
$\lambda \ge
\lambda_*$ and  a.e. $y \in D(\lambda^{-1} 2^{-1} c_4^{-1/2}\delta_0,   \lambda^{-1} \delta_0)$,
\begin{equation}\label{e2.1}
(\Delta+ \sA) g_{\lambda} (y)  \ge \Delta g_{\lambda} (y)- | \sA
g_{\lambda} (y)| \ge \lambda^2 \delta_0^{-2}-4^{-1} \lambda^2 \delta_0^{-2}-c_4
\lambda^4
 \delta_0^{-2} \rho(y)^2
\ge \frac{1}{2} \lambda^2 \delta_0^{-2}
\end{equation}
and
\begin{equation}\label{e:ssp6}
(\Delta+ \sA) h_\lambda(y) \le \lambda|(\Delta+ \sA) h(y)| -
(\Delta+ \sA)g_{\lambda} (y) \le  \lambda (C_1 -\frac{1}{2} \lambda
\delta_0^{-2}) \le - \frac{1}{4} \lambda^2 \delta_0^{-2}\, .
\end{equation}

Let $\delta_*:=2^{-1} c_4^{-1/2}\delta_0$ and
  $f$ be a non-negative smooth radial function with compact
support such that $f(x)=0$ for $|x|>1$ and $\int_{\R^d} f (x) dx=1$.
For $k\geq 1$, define $f_k(x)=2^{kd} f (2^k x)$ and
$$
h_\lambda^{(k)}(z):= ( f_k*h_\lambda)(z) :=\int_{\R^d}f_k (y)
h_\lambda(z-y)dy\, .
$$
Let
$$
B_k^\lambda:=\left\{y \in D(\lambda^{-1} \delta_*, \lambda^{-1}
\delta_0): \delta_{D(\lambda^{-1} \delta_*, \lambda^{-1} \delta_0)
}(y) \ge 2^{-k}\right\}
$$
and consider large $k$'s such that $B_k^\lambda$'s are non-empty
open sets. Since $h_\lambda^{(k)}$ is in $C_c^{\infty}$, $\sA
h_\lambda^{(k)}$ is well defined everywhere. We claim that for every
$\lambda \ge \lambda_*$ and $k$ large enough,
\begin{equation} \label{e:*333}
(\Delta+\sA) h_\lambda^{(k)} \leq -\frac{1}{4} \lambda^2
\delta_0^{-2} \quad \text{ on } B_k^\lambda\, .
\end{equation}
Indeed, for any $x\in B_k^\lambda$ and $z\in B(0,2^{-k})$,  when $k$
is large enough, it holds that $x-z\in D(\lambda^{-1} \delta_*,
\lambda^{-1} \delta_0)$. By
 the proof of Lemma \ref{L:Main} the following limit exists:
\begin{eqnarray*}
\lefteqn{\lim_{\varepsilon\to 0} \int_{B(x,\varepsilon)^c}
\left(h_\lambda(y-z)-h_\lambda(x-z)\right)\, j(|x-y|)\, dy}\\
&=&\lim_{\varepsilon\to 0} \int_{B(x-z, \varepsilon)^c}
\left(h_\lambda(y')-h_\lambda(x-z)\right)\, j(|(x-z)-y'|)\, dy'
\,=\,\sA h_\lambda(x-z)\, .
\end{eqnarray*}
Moreover, by  \eqref{e:ssp6} it holds that for every $\lambda \ge
\lambda_*$,
$(\Delta+ \sA) h_\lambda \le -\frac{1}{4} \lambda^2 \delta_0^{-2}$
a.e.
 on $D(\lambda^{-1} \delta_*, \lambda^{-1} \delta_0)$. Next,
\begin{eqnarray*}
\lefteqn{\int_{B(x, \varepsilon)^c} (h_\lambda^{(k)}(y)-
h_\lambda^{(k)}(x))\, j(|x-y|)\, dy}\\
&=&\int_{|x-y|>\varepsilon} \left(\int_{\R^d}f_k(z)
(h_\lambda(y-z)-h_\lambda(x-z))\, dz\right)\, j(|x-y|)\, dy\\
&=&\int_{B(0, 2^{-k})}f_k(z) \left(\int_{B(x, \varepsilon)^c}
\left(h_\lambda(y-z)-h_\lambda(x-z)\right)\, j(|x-y|)\, dy\right)\,
dz.
\end{eqnarray*}
By letting $\varepsilon \to 0$ and using the dominated convergence
theorem, we get that for every
$\lambda \ge \lambda_*$ and $k$ large enough,
$$
(\Delta + \sA) h_\lambda^{(k)}(x)=\int_{|z|<2^{-k}}
f_k(z)(\Delta+\sA) h_\lambda(x-z)\, dz\le  -\frac{1}{4} \lambda^2
\delta_0^{-2}\int_{|z|<2^{-k}} f_k(z) \, dz =-\frac{1}{4} \lambda^2
\delta_0^{-2}\, .
$$

By using Dynkin's formula \eqref{e:*334}, the estimates \eqref{e:*333} and the fact
that $h_\lambda^{(k)}$ are in $C^\infty_c(\R^d)$, and by letting
$k\to \infty$ we get for every $\lambda \ge
\lambda_*$ and $x \in
D(\lambda^{-1} \delta_*, \lambda^{-1} \delta_0)$ with $\wt x=0$,
\begin{eqnarray}\label{e:ssp7}
\E_x[h_\lambda (X_{\tau_{D(\lambda^{-1} \delta_*, \lambda^{-1}
\delta_0)}})]-\lambda V( \delta_D(x)) &\le & \E_x[h_\lambda
(X_{\tau_{D(\lambda^{-1} \delta_*, \lambda^{-1}
\delta_0)}})]-h_\lambda(x)\nonumber \\
& \le & - \frac{1}{4} \lambda^2 \delta_0^{-2}
\E_x[\tau_{D(\lambda^{-1} \delta_*, \lambda^{-1} \delta_0)}]\, .
 \end{eqnarray}
It is easy to see that  $h_\lambda \ge 0$. In fact, if  $y \in (B(0,
R) \cap D)^c$, then both $h(y)$ and $g_\lambda(y)$ are zero.  If
$y \in B(0, R) \cap D$ and $\rho(y ) \ge 4 \lambda^{-1} \delta_0$,
then $g_\lambda(y)=0$. Finally, if $y \in B(0, R) \cap D$ and
$\rho(y ) \le 4 \lambda^{-1} \delta_0$, then, since $g(y) \le y_d^2$
by \eqref{e:ssp101}, we have from \eqref{e:ssp1},
$$
h_\lambda(y) = \lambda V(\delta_D(y)) - g( \lambda \delta_0^{-1} \wt
y,  \lambda \delta_0^{-1} \rho(y)) \ge \lambda V(\delta_D(y)) -
\lambda^2 \rho(y)^2 \ge 0.
$$
Therefore, from \eqref{e:ssp7},
 \bee\label{e:ssp8}
V( \delta_D(x)) \,\ge\,  \frac{1}{4}\, \lambda \,\delta_0^{-2}
\,\E_x[\tau_{D(\lambda^{-1} \delta_*, \lambda^{-1} \delta_0)}].
 \eee
Since $B(0, (1+\Lambda)^{-1} \delta_* \lambda^{-1}) \cap D \subset
D(\lambda^{-1} \delta_*, \lambda^{-1} \delta_0)$, using Lemma
\ref{l2.1} and \eqref{e:ssp8}, we have that for every $\lambda \ge
\lambda_*$ and $x \in B(0, 2^{-1}(1+\Lambda)^{-1} \delta_*
\lambda^{-1})$ with $
\wt x=0$,
$$
\P_{x}\left(X_{ \tau_{ D_Q (  \lambda^{-1} \delta_*,  \lambda^{-1}
\delta_0)}} \in  D\right) \,\le\, c_7\,  \lambda^2
\,\E_x[\tau_{D(\lambda^{-1} \delta_*, \lambda^{-1} \delta_0)}] \,\le
\,c_8 \,\lambda\, V( \delta_D(x)).
$$
We have proved the lemma with $\lambda_0:=\lambda_*\delta_0^{-1}$.
\qed

\begin{lemma}\label{L:200}
There is a constant $c=c(R,  \Lambda )>0$ such that for every
$\lambda \ge \lambda_0$, $\kappa \in (0, 1]$,  $Q \in
\partial D$ and $x \in D_Q ( \kappa\lambda^{-1} ,  \lambda^{-1} )$ with $\wt x=0$,
 \bee\label{e:L:1}
\P_{x}\left(X_{ \tau_{ D_Q ( \kappa \lambda^{-1} ,  \lambda^{-1} )}} \in
 D_Q ( 2 \kappa\lambda^{-1} ,  \lambda^{-1} )\right) \ge c  \lambda V(
\delta_D (x)).
 \eee
\end{lemma}

\pf  Fix $\lambda \ge \lambda_0$ and $\kappa \in (0, 1]$.
For simplicity we denote $D_Q(\kappa\lambda^{-1} ,  \lambda^{-1})$ by $\wh{D}$.
Further, let
$$
B=\left\{y\in D: \rho_Q(y)=\kappa\lambda^{-1}
\mbox{ and } |\wt{y}|<
\lambda^{-1}\right\}
$$
be the upper boundary of $\wh{D}$.

Let $\tau^W_{\wh{D}}$ be the first time the Brownian motion $W$ exits $\wh{D}$
and $W^{\wh{D}}$
be the killed Brownian motion in $\wh{D}$.
Let $Y=(Y_t:\, t\ge 0)$ be the subordinate killed Brownian motion defined by
$Y_t=W^{\wh{D}}_{S_t}$.
Let $\zeta$ denote the
lifetime of $Y$.
Recall that $u$ is the potential density of the subordinator $S$. It
follows from \cite[Corollary 4.4]{SV08} that
$$
\P_x(X_{\tau_{\wh{D}}}\in B)\,
\ge\, \P_x(Y_{\zeta-}\in B)\,=\,\E_x\left[u(\tau^W_{\wh{D}});
W_{\tau^W_{\wh{D}}}\in B\right].
$$
Thus, since $u$ is deceasing, for any $t>0$,
\begin{align*}
&\P_x(X_{\tau_{\wh{D}}}\in B)\,\ge \,\E_x\left[u(\tau^W_{\wh{D}});
W_{\tau^W_{\wh{D}}}\in B, \tau^W_{\wh{D}}\le t \right] \,\ge\,
u(t)\P_x\big(W_{\tau^W_{\wh{D}}}\in B, \tau^W_{\wh{D}}\le t \big)\\
&=\, u(t)\left[\P_x\big(W_{\tau^W_{\wh{D}}}\in B\big)-\P_x\big(\tau^W_{\wh{D}}> t\big)\right] \,\ge \,
u(t)\left[\P_x\big(W_{\tau^W_{\wh{D}}}\in B\big) -t^{-1}
\E_x\big[\tau^W_{\wh{D}}\big]\right]\, .
\end{align*}

Now we use the following two estimates which
are valid for the Brownian motion (for example, see \cite[Lemma
3.4]{CKSV1} with $a=0$). There exist constants $c_1>0$ and
$c_2>0$ (independent of $\lambda \ge  \lambda_0 $) such that
$
\P_x\big(W_{\tau^W_{\wh{D}}}\in B \big) \ge c_1 \lambda
\delta_D(x)$ and $ \E_x\big[\tau^W_{\wh{D}}
\big]\le c_2 \lambda^{-1} \delta_D(x)\, .
$
Then, by choosing $t_0>0$ so that $c_1-t^{-1}_0 c_2 \lambda^{-2} \ge
c_1-t^{-1}_0 c_2 \lambda_0  \ge c_1/2=:c_3$, we get
$$
\P_x\big(X_{\tau_{\wh{D}}}\in B \big)\ge u(t)(c_1   -c_2
t^{-1}\lambda^{-2})\lambda \delta_D(x)\ge c_3 u(t_0) \lambda
\delta_D(x)\, .
$$\qed

\section{Carleson estimate and Boundary Harnack principle}

In this section, we give the proof of the boundary Harnack principle
for $X$. We first prove the Carleson estimate for $X$ on Lipschitz
open sets.

We recall that an open set $D$ in $\bR^d$ is said to be a Lipschitz
open set if there exist a localization radius $R_2>0$ and a constant
$\Lambda_1 >0$ such that for every $Q\in
\partial D$, there exist a Lipschitz function $\psi=\psi_Q:
\bR^{d-1}\to \bR$ satisfying $\psi (0)= 0$, $| \psi (x)- \psi (y)|
\leq \Lambda_1 |x-y|$, and an orthonormal coordinate system $CS_Q$:
$y=(y_1, \dots, y_{d-1}, y_d)=:(\wt y, \, y_d)$ with its origin at
$Q$ such that
$$
B(Q, R_2)\cap D=\{ y=(\wt y, y_d)\in B(0, R_2) \mbox{ in } CS_Q: y_d
> \psi (\wt y) \}.
$$
The pair $(R_2, \Lambda_1)$ is called the characteristics of the
Lipschitz open set $D$.
Without loss of generality, we will assume throughout this section that
$R_1<1$.
Note that a Lipschitz open set can be
unbounded and disconnected. For Lipschitz open set $D$ and every
$Q\in \partial D$ and $ x \in B(Q, R_2)\cap D$, we define
$$
\rho_Q (x) := x_d -  \psi_Q (\tilde x)\, ,
$$
where $(\tilde x, x_d)$ are the coordinates of $x$ in $CS_Q$.

The proof of the next lemma is similar to that of \cite[Lemma 4.1]{CKSV1}.
\begin{lemma}\label{lower bound} Let $D\subset \R^d$ be a Lipschitz
open set with characteristics $(R_2, \Lambda_1)$.  There exists
a constant $\delta=\delta(R_2, \Lambda_1)>0$ such that for all
$Q \in \partial D$ and  $x\in D$ with $\rho_Q(x) < R_2/2$,
$$
\P_x(X_{\tau(x)}\in D^c)\ge \delta\, ,
$$
where $\tau(x):=\tau_{D\cap B(x,2\rho_Q(x))}=\inf\{t>0:\,
X_t\notin D\cap B(x,2\rho_Q(x))\}$.
\end{lemma}

\pf
Let $D_x:=D \cap B(x,2\rho_Q(x))$ and $W^{D_x}$ be the killed Brownian
motion in $D_x$. Here $W$ denotes the Brownian motion in $\R^d$. As in the proof of
Lemma \ref{L:200}, we define the subordinate killed Brownian motion
$Y=(Y_t:\, t\ge 0)$ in $D_x$ by $Y_t:=W^{D_x}(S_t)$.
We will use $\zeta$ to denote the
lifetime of $Y$ and let $C_x:=\partial D\cap
B(x,2\rho_Q(x))$
and $\tau^W _{U} :=\inf\{t>0:\, W_t\notin {U}\}$.

Since, see \cite{SV08},
$
\P_x\left(X_{\tau(x)}\in C_x\right)\ge
\P_x\left(Y_{\zeta-}\in C_x\right)=\E_x\left[u(\tau^W _{D_x});\, W_{\tau^W _{D_x}}\in
C_x\right],
$
we have
\begin{eqnarray}\label{ineq1}
&&\P_x\left(X_{\tau(x)}\in D^c\right)\, \ge\, \P_x\left(X_{\tau(x)}\in
C_x\right) \ge
\E_x\left[u(\tau^W _{D_x}); \, W_{\tau^W _{D_x}}\in C_x,
\tau^W _{D_x} \le t\right]\nonumber \\
&&\ge  u(t)\P_x\left[W_{\tau^W _{D_x}}\in C_x, \tau^W _{D_x}
\le t\right] \ge  u(t)\left(\P_x(W_{\tau^W _{D_x}}\in
C_x)-\P_x(\tau^W _{D_x}>t)\right), \quad t>0.
\end{eqnarray}
By the fact that $D$ is a Lipschitz open set,  there exists
$c_1=c_1(R_2, \Lambda_1)>0$ such that
\begin{equation}\label{ineq2}
\P_x(W_{\tau^W _{D_x}}\in C_x)\ge c_1\, .
\end{equation}
(See the proof of \cite[Lemma 4.1]{CKSV1}.)
Since
$$
\P_x(\tau^W _{D_x} >t)\le \frac{\E_x[\tau^W _{D_x}]}{t} \le
\frac{\E_x[\tau^W _{B(x, 2\rho_Q(x))}]}{t} \le c_2
\frac{(\rho_Q(x))^2}{t}\le c_2\frac{R_2^2}{t},$$
by using (\ref{ineq2}) and (\ref{ineq1}), we obtain that
$$
\P_x\left(X_{\tau(x)}\in D^c\right) \ge  u(t)\left(\P_x(W_{\tau^W _{D_x}}\in
C_x)-\P_x(\tau^W _{D_x}>t)\right)
 \ge u(t)  \left(c_1-c_2\frac{R^2_1}{t}\right) \ge  c_1 u(t_0)/2 >0,
$$
where $t_0=t_0(R_2, \Lambda_1)>0$ is chosen so that $c_1-c_2
R_2^2/t\ge c_1/2$. The lemma is thus proved. \qed

Suppose that $D$ is an open set and that $U$ and $V$ are bounded
open sets with $V \subset \overline{V} \subset U$ and $ D \cap V
\not= \emptyset$. If $f$ vanishes continuously on $D^c\cap U$, then
by a finite covering argument, it is easy to see that $f$ is bounded
in an open neighborhood of $\partial D\cap V$. The proof of the next
result is the same as that of \cite[Lemma 4.2]{CKSV1}. So we omit
the proof.

\begin{lemma}\label{l:regularity}
Let $D$ be an open set and $U$ and $V$ be bounded open sets with $V
\subset \overline{V} \subset U$ and $ D \cap V \not= \emptyset$.
Suppose $f$ is a nonnegative function in $\R^d$ that is harmonic in
$D\cap U$ with respect to $X$ and vanishes continuously on $D^c\cap
U$. Then $f$ is regular harmonic in $D\cap V$ with respect to $X$,
i.e.,
\begin{equation}\label{e:regularity}
f(x)=\E_x\left[ f(X_{\tau_{D\cap V}})\right] \qquad \hbox{
for all }x\in D\cap V\, .
\end{equation}
\end{lemma}

\begin{thm}[Carleson estimate]\label{carleson}
Let $D\subset \R^d$ be a Lipschitz open set with the characteristics
$(R_2, \Lambda_1)$. Then there exists a positive constant
$A=A(R_2, \Lambda_1)$ such that for every $Q\in
\partial D$, $0<r<R_2/2$, and any nonnegative function
$f$ in $\R^d$ that is harmonic in $D \cap B(Q, r)$ with respect to
$X$ and vanishes continuously on $ D^c \cap B(Q, r)$, we have
\begin{equation}\label{e:carleson}
f(x)\le A f(x_0) \qquad \hbox{for }  x\in D\cap B(Q,r/2),
\end{equation}
where $x_0\in D
\cap B(Q,r)$ with $\rho_Q(x_0)=r/2$.
\end{thm}

\pf Since $D$ is Lipschitz and $r<R_2/2$, by
Proposition \ref{uhp}
and a standard chain argument, it suffices to prove
(\ref{e:carleson}) for $x\in D\cap B(Q,r/12)$ and $\wt x_0 = \wt Q$.
Without loss of generality, we may assume that $f(x_0)=1$. In this
proof, the constants $\delta, \beta, \eta$ and $c_i$'s are always
independent of $r$.

Let $\nu=\nu(3) \vee 2 $ where $\nu(3)$ is the constant in \eqref{H:1n} with $K=3$,
choose $0<\gamma < (\nu^{-1} \wedge (1-\nu^{-1}))$ and let
$$
B_0(x)=D\cap B(x,2\rho_Q(x))\, ,\qquad B_1(x)=B(x,r^{1-\gamma}
\rho_Q(x)^{\gamma})\,
$$
and
$$
B_2=B(x_0,\rho_Q(x_0)/3)\, ,\qquad B_3=B(x_0, 2\rho_Q(x_0)/3).
$$
By Lemma \ref{lower bound}, there exists $\delta=\delta(R_2,
\Lambda_1)>0$ such that
\begin{equation}\label{e:c:1}
\P_x(X_{\tau_{B_0(x)}}\in D^c)\ge \delta\, ,\quad x\in B(Q,r/4)\, .
\end{equation}
By the Harnack inequality and a chain argument, there exists
$\beta>0$ such that
\begin{equation}\label{e:c:2}
f(x)<(\rho_Q(x)/r)^{-\beta} f(x_0)\, ,\quad x\in D\cap B(Q,r/4)\, .
\end{equation}
In view of Lemma \ref{l:regularity}, $f$ is regular harmonic in
$B_0(x)$ with respect to $X$. So
\begin{equation}\label{e:c:3}
f(x)=\E_x\big[f\big(X_{\tau_{B_0(x)}}\big); X_{\tau_{B_0(x)}}\in
B_1(x)\big]+ \E_x\big[f\big(X_{\tau_{B_0(x)}}\big);
X_{\tau_{B_0(x)}}\notin B_1(x)\big] \qquad \hbox{for } x\in B(Q,
r/4) .
\end{equation}
We first show that there exists $\eta>0$ such that
\begin{equation}\label{e:c:4}
\E_x\big[f\big(X_{\tau_{B_0(x)}}\big); X_{\tau_{B_0(x)}}\notin
B_1(x)\big]\le f(x_0)  \quad \hbox{if } x\in D \cap B(Q,r/12) \hbox{
with } \rho_Q(x) < \eta r\, .
\end{equation}
Let $\eta_0 :=2^{-2 \nu }$,
then,
since $\gamma < 1-\nu^{-1}$,
for $\rho_Q(x)< \eta_0 r$,
$$
2\rho_Q(x) \le r^{1-\gamma} \rho_Q(x)^{\gamma} - 2\rho_Q(x).
$$
Thus if $x\in D \cap B(Q,r/12)$ with $\rho_Q(x) < \eta_0r$,  then
$|x-y|\le 2|z-y|$ for $z\in B_0(x)$, $y\notin B_1(x)$. Moreover, by
the triangle inequality, $|x-y|\le |x-z|+|z-y|\le 1+|z-y|$. Thus we
have by \eqref{H:1}, \eqref{H:2}, \eqref{e:levy} and Lemma
\ref{L:2.00}
\begin{align}\label{e:c:5}
&\E_x\big[f\big(X_{\tau_{B_0(x)}}\big); X_{\tau_{B_0(x)}} \notin
B_1(x)\big]\nonumber\\
=&\E_x \int_0^{\tau_{B_0(x)}}   \int_{2>|y-x|>r^{1-\gamma}
\rho_Q(x)^{\gamma}}j(|X_t-y|)f(y)\, dy\, dt
\nonumber +\E_x \int_0^{\tau_{B_0(x)}}\int_{|y-x|>2}j(|X_t-y|)f(y)\, dy\, dt
\nonumber \\
\le &c_1 \E_x [\tau_{B_0(x)}]\left(\int_{2>|y-x|>r^{1-\gamma}
\rho_Q(x)^{\gamma}}j(|x-y|)f(y)\, dy+\int_{|y-x|>2}j(|x-y|)f(y)\, dy\right)
\nonumber \\
\le &c_1  c_2 \rho_Q(x)^2
\left(\int_{|y-x|>r^{1-\gamma}\rho_Q(x)^{\gamma}, |y-x_0|>2
\rho_Q(x_0)/3}
j(|x-y|)f(y)\, dy \right.\nonumber \\
&\quad \quad \quad \quad \quad +\left.\int_{|y-x_0|\le
2\rho_Q(x_0)/3}j(|x-y|)f(y)\, dy\right)\,=:\, c_3 \rho_Q(x)^2
(I_1+I_2)\, .
\end{align}
On the other hand, for $z\in B_2$ and $y\notin B_3$, we have
$|z-y|\le |z-x_0|+|x_0-y|\le \rho_Q(x_0)/3+|x_0-y|\le 2|x_0-y|$ and
$|z-y|\le |z-x_0|+|x_0-y|\le 1+|x_0-y|$. We have again by
\eqref{e:levy}, \eqref{H:1}, \eqref{H:2} and Lemma \ref{L:2.00}
\begin{align}\label{e:c:6}
&f(x_0)\,\ge\, \E_{x_0}\left[f(X_{\tau_{B_2}}), X_{\tau_{B_2}}\notin B_3\right]\nonumber \\
&\ge    \E_{x_0} \int_0^{\tau_{B_2}}  \left(\int_{2>|y-x_0|>2\rho_Q(x_0)/3} j(|X_t-y|)f(y)\, dy
+ \int_{|y-x_0|\ge 2} j(|X_t-y|)f(y)\, dy\right) dt\nonumber \\
&\ge c_4\E_{x_0} [\tau_{B_2}]    \left(\int_{2>|y-x_0|>2\rho_Q(x_0)/3} j(|x_0-y|)f(y)\, dy
+ \int_{|y-x_0|\ge 2} j(|x_0-y|)f(y)\, dy\right)\nonumber \\
&=c_5 \rho_Q(x_0)^2 \int_{|y-x_0|>2\rho_Q(x_0)/3} j(|x_0-y|)f(y)\, dy\, .
\end{align}

Suppose now that $|y-x|\ge r^{1-\gamma}\rho_Q(x)^{\gamma}$ and $x\in
B(Q,r/4)$. Then
$$
|y-x_0|\le |y-x|+r\le
|y-x|+r^{\gamma}\rho_Q(x)^{-\gamma}|y-x|\le 2r^{\gamma}\rho_Q(x)^{-\gamma}|y-x|.
$$
Thus, using \eqref{H:1n}, we
get  for $|x-y| \le 2$,
\begin{eqnarray}\label{e:gf1} j(|y-x| )
 \le  c_7 (\rho_Q(x)/r)^{-\nu \gamma} j(|y-x_0| ).
\end{eqnarray}
Now, using \eqref{H:1}, \eqref{H:2} and \eqref{e:gf1},
\begin{eqnarray}\label{e:c:7}
I_1&\le &c_7 \int_{R_0/2>|y-x|>r^{1-\gamma}\rho_Q(x)^{\gamma},
|y-x_0|>2\rho_Q(x_0)/3}(\rho_Q(x)/r)^{-
\nu \gamma}
j(|y-x_0| ) f(y)\, dy\nonumber\\
& &+c_8\int_{|y-x|\ge R_0/2, |y-x_0|>2\rho_Q(x_0)/3}j(|x_0-y|)f(y)\,
dy\nonumber\\
&\le & c_9  \left( (\rho_Q(x)/r)^{-
\nu \gamma}+1\right)
\int_{ |y-x_0|>2\rho_Q(x_0)/3}j(|x_0-y|)\, f(y)\, dy\nonumber\\
&\le &c_5^{-1} c_9\rho_Q(x_0)^{-2}\left( (\rho_Q(x)/r)^{-
\nu \gamma}+1\right)f(x_0)\nonumber\\
&\le& 2c_5^{-1} c_9 (\rho_Q(x)/r)^{-
\nu \gamma}
\rho_Q(x_0)^{-2} f(x_0) \, ,
\end{eqnarray}
where the second to last inequality is due to \eqref{e:c:6}.

If $|y-x_0|<2\rho_Q(x_0)/3$, then $|y-x|\ge |x_0-Q|-|x-Q|-|y-x_0|
>\rho_Q(x_0)/6$. This together with the Harnack inequality implies that
\begin{eqnarray}\label{e:c:8}
I_2 &\le& c_{10} \int_{|y-x_0|\le 2\rho_Q(x_0)/3}j(|x-y|) f(x_0)\, dy
\le  c_{10} f(x_0)\int_{|y-x|>\rho_Q(x_0)/6} j(|x-y|)\, dy\nonumber\\
&=&  c_{10} f(x_0)\left(\int_{R_0>|z|>\rho_Q(x_0)/6} j(|z|)\, dz +
\int_{R_0 \le |z|} j(|z|)\, dz \right)\nonumber\\
 &\le &  c_{10} f(x_0)\left(\int_{R_0>|z|>\rho_Q(x_0)/6} j(|z|)\, dz +
c_{11} \right).
\end{eqnarray}

Combining \eqref{e:c:5}, \eqref{e:c:7} and \eqref{e:c:8} we obtain
\begin{eqnarray}\label{e:c:9}
&&\E_x[f(X_{\tau_{B_0(x)}});\, X_{\tau_{B_0(x)}}\notin B_1(x)]\nonumber\\
&\le &c_{12} f(x_0)\Big(\rho_Q(x)^2(\rho_Q(x)/r)^{-
\gamma
\nu}\rho_Q(x_0)^{-2}
\nonumber \\
&&\quad
(\rho_Q(x)/r)^2(\rho_Q(x_0)/6)^2  \int_{R_0>|z|>\rho_Q(x_0)/6} j(|z|)\, dz     +(\rho_Q(x)/r)^2 r^2 \Big)\nonumber \\
&\le &c_{13} f(x_0)\left((\rho_Q(x)/r)^{2-
\gamma
\nu}
+  (\rho_Q(x)/r)^2 \Big(\int_{R_0>|z|>\rho_Q(x_0)/6} |z|^2j(|z|)\, dz     +1\Big) \right)\nonumber \\
&\le &c_{14} f(x_0)\left((\rho_Q(x)/r)^{2-
\gamma
\nu}
  +  (\rho_Q(x)/r)^2   \right),
\end{eqnarray}
where  we used the fact that
$\rho_Q(x_0)=r/2$.  Since
$2-\gamma
\nu>0$, choose now $\eta\in (0, \eta_0)$ so that
$$
c_{14}\,\left(\eta^{2-\gamma
\nu} +\eta^2
\right)\,\le\, 1\, .
$$
Then  for $x\in  D \cap B(Q,r/12)$ with $\rho_Q(x) < \eta r$, we
have by \eqref{e:c:9},
\begin{eqnarray*}
\E_x\left[f(X_{\tau_{B_0(x)}});\, X_{\tau_{B_0(x)}}\notin
B_1(x)\right] &\le & c_{14}\,
f(x_0)\left(\eta^{2-\gamma
\nu}+\eta^2 \right) \le
f(x_0)\, .
\end{eqnarray*}

We now prove the Carleson estimate \eqref{e:carleson} for $x\in
D\cap B(Q, r/12)$ by a method of contradiction. Recall that
$f(x_0)=1$. Suppose that there exists $x_1\in D\cap B(Q,r/12)$ such
that $f(x_1)\ge K>\eta^{-\beta}\vee (1+\delta^{-1})$, where $K$ is a
constant to be specified later. By \eqref{e:c:2} and the assumption
$f(x_1)\ge K>\eta^{-\beta}$, we have
$(\rho_Q(x_1)/r)^{-\beta}>f(x_1)\ge K> \eta^{-\beta}$, and hence
$\rho_Q(x_1)<\eta r$.
By (\ref{e:c:3}) and  (\ref{e:c:4}),
$$
K\le f(x_1)\le \E_{x_1}\left[f(X_{\tau_{B_0(x_1)}});
X_{\tau_{B_0(x_1)}} \in B_1(x_1) \right]+1\, ,
$$
and hence
$$
\E_{x_1}\left[f(X_{\tau_{B_0(x_1)}}); X_{\tau_{B_0(x_1)}} \in
B_1(x_1)\right] \ge f(x_1)-1 > \frac{1}{1+\delta}\, f(x_1)\, .
$$
In the last inequality of the display above we used the assumption
that  $f(x_1)\ge K>1+\delta^{-1}$. If $K \ge 2^{\beta/\gamma}$, then
$D^c\cap B_1(x_1)\subset D^c \cap B(Q,r)$. By using the assumption
that $f=0$ on $D^c\cap B(Q, r)$, we get from \eqref{e:c:1}
\begin{eqnarray*}
\E_{x_1}[f(X_{\tau_{B_0(x_1)}}), X_{\tau_{B_0(x_1)}}\in B_1(x_1)]
&=&\E_{x_1}[
f(X_{\tau_{B_0(x_1)}}), X_{\tau_{B_0(x_1)}}\in B_1(x_1)\cap D]\\
& \le& \P_x(X_{\tau_{B_0(x_1)}}\in
D) \, \sup_{B_1(x_1)}f  \le (1-\delta) \, \sup_{B_1(x_1)}f \, .
\end{eqnarray*}
Therefore, $\sup_{B_1(x_1)}f> f(x_1)/((1+\delta)(1-\delta))$, i.e.,
there exists a point $x_2\in D$ such that
$$
|x_1-x_2|\le r^{1-\gamma}\rho_Q(x_1)^{\gamma} \quad \hbox{ and }
\quad
f(x_2)>\frac{1}{1-\delta^2}\, f(x_1)\ge \frac{1}{1-\delta^2}\, K\, .
$$
By induction, if $x_k\in D\cap B(Q, r/12)$ with $f(x_k)\geq
K/(1-\delta^2)^{k-1}$ for $k\ge 2$, then there exists $x_{k+1}\in D$
such that
\begin{equation}\label{e:c:10}
|x_k-x_{k+1}|\le r^{1-\gamma}\rho_Q(x_k)^{\gamma}  \quad \hbox{ and
} \quad f(x_{k+1}) > \frac{1}{1-\delta^2}\, f(x_k)>
\frac{1}{(1-\delta^2)^k}\, K\, .
\end{equation}
{}From (\ref{e:c:2}) and (\ref{e:c:10}) it follows that
$\rho_Q(x_{k})/r \le (1-\delta^2)^{(k-1)/\beta}K^{-1/\beta}$, for
every $k\ge 1$. Therefore,
\begin{align*}
&|x_k-Q|\,\le\,|x_1-Q|
+\sum_{j=1}^{k-1}|x_{j+1}-x_j|\,\le\, \frac{r}{12} +
 \sum_{j=1}^{\infty} r^{1-\gamma}\rho_Q(x_j)^{\gamma}\\
&\le \frac{r}{12}+r^{1-\gamma}\sum_{j=1}^{\infty}(1-\delta^2)^{(j-1)
\gamma/\beta}K^{-\gamma/\beta}r^{\gamma}\,=\,\frac{r}{12}+
r^{1-\gamma}r^{\gamma}K^{-\gamma/\beta}
\sum_{j=0}^{\infty}(1-\delta^2)^{j\gamma/\beta}\\
&=\frac{r}{12}+ r K^{-\gamma/\beta}\,
\frac{1}{1-(1-\delta^2)^{\gamma/\beta}}.
\end{align*}
Choose
$$
K=\eta\vee (1+\delta^{-1})\vee 12^{\beta/\gamma}(1-
(1-\delta^2)^{\gamma/\beta})^{-\beta/\gamma}.
$$
Then $K^{-\gamma/\beta}\, (1-(1-\delta^2)^{\gamma/\beta})^{-1}\le
1/12$, and hence $x_k\in D\cap B(Q,r/6)$ for every $k\ge 1$. Since
$\lim_{k\to \infty}f(x_k)=+\infty$, this contradicts the fact that
$f$ is bounded on $B(Q,r/2)$. This contradiction shows that $f(x)<
K$ for every $x\in D\cap B(Q, r/12)$. This completes the proof of
the theorem.
 \qed

\noindent {\bf Proof of Theorem \ref{t:main} }.
 We recall  that $R_1=R/(4\sqrt{1 + (1+\Lambda)^2})$ and $\lambda_0
>2 R_1^{-1}$ and $\kappa_0 \in (0,1)$ are the constants in the statement of
Lemma \ref{L:2}.

Since $D$ is a $C^{1,1}$ open set and $r<R$, by the Harnack
inequality and a standard chain argument, it
suffices to prove \eqref{e:bhp_m} for $x,y \in D \cap B(Q,2^{-1}
r\kappa_0\lambda_0^{-1})$. In this proof, the constants $\eta$ and $c_i$'s
are always independent of $r$.

For any $r\in (0, R]$ and
 $x\in D\cap B(Q, 2^{-1} r\kappa_0\lambda_0^{-1})$, let $Q_x$ be
the point $Q_x \in \partial D$ so that $|x-Q_x|=\delta_{D}(x)$ and
let $x_0:=Q_x+\frac{r}{8}(x-Q_x)/|x-Q_x|$. We choose a
$C^{1,1}$-function $\varphi: \bR^{d-1}\to \bR$ satisfying $\varphi
(0)= 0$, $\nabla\varphi (0)=(0, \dots, 0)$, $\| \nabla \varphi
 \|_\infty \leq \Lambda$, $| \nabla \varphi (y)-\nabla \varphi (z)|
\leq \Lambda |y-z|$, and an orthonormal coordinate system $CS$ with
its origin at $Q_x$ such that
$$
B(Q_x, R)\cap D=\{ y=(\wt y, y_d) \in B(0, R) \mbox{ in } CS: y_d >
\varphi (\wt y) \}.
$$
In the coordinate system $CS$ we have $\wt x = \wt 0$ and $x_0=(\wt
0, r/8)$. For any $b_1, b_2>0$, we define
$$
D(b_1, b_2):=\left\{ y=(\wt y, y_d) \mbox{ in } CS: 0<y_d-\varphi(\wt
y)<b_1r\kappa_0\lambda_0^{-1}, \ |\wt y| < b_2 r\lambda_0^{-1} \right\}.
$$
It is easy to see that
$D(2, 2)\subset D\cap B(Q, r/2)$.
In fact, since $\Lambda \ge 1$ and  $R \le 1$,
 for every $z \in  D(2, 2)$,
\begin{align*}
&|z-Q| \le |Q-x|+|x-Q_x|+|Q_x-z| \le |Q-x|+|x-Q_x|+ |z_d- \varphi(\wt
z)|+ |\varphi(\wt z)|\\
&< r \lambda_0^{-1} ((1+\Lambda)+ 4)
<2^{-1} r R((1+\Lambda)+ 4)/(4\sqrt{1 + (1+\Lambda)^2})
  \le \frac{r}{2}.
\end{align*}
Thus if $f$ is a nonnegative function on $\R^d$ that is harmonic in
$D\cap B(Q, r)$ with respect to $X$ and vanishes continuously in
$D^c\cap B(Q, r)$, then, by Lemma \ref{l:regularity}, $f$ is regular
harmonic in $D\cap B(Q,r/2)$ with respect to $X$, hence also in
$D(2, 2)$. Thus by the  Harnack
inequality, we have
\begin{eqnarray}
f(x) &= & \E_x\left[f\big(X_{\tau_{ D(1,1)}}\big)\right] \ge
\E_x\left[f\big(X_{\tau_{ D(1,1)}}\big); X_{
\tau_{ D(1,1)}} \in  D(2,1)\right]\label{e:BHP2}\\
&\ge& c_1 f(x_0) \P_x\Big( X_{\tau_{ D(1,1)}} \in D (2,1)\Big)
\ge c_2 f(x_0) \delta_D(x) /r.\nonumber
\end{eqnarray}
In the last inequality above we have used \eqref{e:L:1}.

Let $w=(\wt 0, r\lambda_0^{-1}\kappa_04)$. Then it is easy to see that there
exists a constant $\eta=\eta(\Lambda,  \delta_0)\in (0, 1/4)$ such
that $B(w, \eta r\lambda_0^{-1}\kappa_0)\in D(1, 1)$.  By  \eqref{H:1},
\eqref{H:2}, \eqref{e:levy} and Lemma \ref{L:2.00},
\begin{align*}
&f(w) \,\ge\, \E_{w}\left[f\big(X_{\tau_{ D(1,1)}}\big);
X_{\tau_{ D(1,1)}} \notin  D(2,2)\right]\,=\,
\E_{w}
\int_0^{\tau_{ D(1,1)}}
\int_{\R^d\setminus D(2,2)} f(y)     j(|X_t-y|)dydt\\
&\ge\, c_3\,\E_{w}\big[\tau_{B(w, \eta r\lambda_0^{-1}\kappa_0 )}\big]
\int_{\R^d\setminus  D(2,2)} f(y) j(|w-y|)  dy\,\ge\, c_4\,  r^2\,
\int_{\R^d\setminus  D(2,2)} f(y) j(|w-y|)   dy.
\end{align*}
Hence by  \eqref{H:1}, \eqref{H:2}, \eqref{e:L:3},
\begin{align*}
&\E_{x}\left[f \left(X_{\tau_{ D(1,1)}}\right); \,
X_{\tau_{ D(1,1)}}
\notin  D(2,2)\right]\,=\, \E_{x}
\int_0^{\tau_{ D(1,1)}}
\int_{\R^d\setminus
D(2,2)}  f(y)     j(|X_t-y|)dydt\\
&\le\, c_5\, \E_x[\tau_{ D(1,1)}]  \int_{\R^d\setminus
D(2,2)}   f(y)     j(|w-y|)dy\\
&\le\, c_6\,  \delta_D(x) r  \int_{\R^d\setminus  D(2,2)}
f(y) j(|w-y|)  dy \,\leq \,  \frac{c_6 \,
\delta_D(x)}{c_4 \, r} f(w).
\end{align*}

On the other hand, by the Harnack inequality and the Carleson estimate, we have
$$
\E_x\left[f\left(X_{\tau_{ D(1,1)}}\right);\,  X_{\tau_{
D(1,1)}} \in D(2,2)\right] \,\le\, c_7 \, f(x_0) \P_x\left(
X_{\tau_{D (1,1)}} \in  D(2,2)\right)\,\le\, c_8 \, f(x_0)
\delta_D(x) /r.
$$
In the last inequality above we have used \eqref{e:L:2}.
Combining the two inequalities above, we get
\begin{align}
& f(x) \,= \, \E_x\left[f \left(X_{\tau_{ D(1,1)}}\right); \,
X_{\tau_{ D(1,1)}} \in D(2,2)\right]\,+\,\E_x\left[ f
\left(X_{\tau_{ D(1,1)}}\right); \,
X_{\tau_{ D(1,1)}} \notin  D(2,2)\right] \label{e:BHP1} \\
&\le \, \frac{c_8}{r} \delta_D(x) f(x_0)   + \frac{c_6 \,
\delta_D(x)}{c_4\, r} f(w) \,\le  \, \frac{c_{9}}{r}\, \delta_D(x)
(f(x_0)  + f(w))\,\le \, \frac{c_{10}}{r}\, \delta_D(x) f(x_0)  .\nonumber
\end{align}
In the last inequality above we have used the  Harnack inequality.

From \eqref{e:BHP2}--\eqref{e:BHP1}, we have that for every $x, y\in
  D \cap B(Q, 2^{-1} r\kappa_0\lambda_0^{-1})$,
$$
\frac{f(x)}{f(y)}\,\le \,
\frac{c_{10}}{c_2}\,\frac{\delta_D(x)}{\delta_D(y)},
$$
which  proves the theorem. \qed

\section{Counterexample}\label{counterexample}

In this section,
we present an example of a (bounded) $C^{1,1}$ domain (open and connected) $D$
on which the boundary Harnack principle for the independent sum of a Brownian motion and
a finite range rotationally invariant L\'evy process fails,
even for regular harmonic function vanishing on $D^c$.
A similar example appears in \cite[Section 6]{KS} for the case of truncated stable process.

Suppose that $Z$ is a rotationally invariant L\'evy process whose L\'evy measure
has a density $J(x)=j(|x|)$ with $j(r)=0$ for all $r\ge 1$. Suppose that $Z$ is independent
of the Brownian motion $W$. We will consider the process $Y=W+Z$.
For any Borel sets $U$ and $V$ in $\R^d$ with $V \subset \overline{U}^c$, we have
\begin{equation}
\label{Poisson}
\P_x(Y_{\tau^Y_U} \in V)=\E_x \int_0^{\tau_U^Y}
\int_{V } j(|Y_t-z|){\bf 1}_{\{|Y_t-z| <1\}}(|Y_t-z|) \, dz\, dt
\quad x \in U,
\end{equation}
where $\tau_U^Y:=\inf\{t>0: Y_t \notin U\}$.

We consider the bounded domain in $\R^d$
$$
D:=(-100, 100)^{d} \setminus \left(    (-100, 49]^{d-1} \times[-1/2, 0]\right).
$$

Suppose that the (not necessarily scale invariant) boundary Harnack principle
is true for $Y$ on $D$ at the origin for regular harmonic function vanishing on $D^c$,
i.e., there exist constants $R_1 >0 $ and  $M_1 >1$ such that for any  $r < R_1$ and
any nonnegative  functions $u, v$ on $\R^d$ which are regular  harmonic with
respect to $Y$ in $D \cap B(
0, M_1 r)$ and vanish in $D^c$, we have
\begin{equation}\label{ce0}
\frac{ u(x)}{ v(x)}
\, \le \,c_r\, \frac{ u(y)}{ v(y)}  \quad  \mbox{ for any }
 x,y \in
D  \cap B(
0, r),
\end{equation}
where $c_r=c_r(D)>0$ is independent of the harmonic functions $u$ and $v$.
Choose an $r_1 < R_1$ with $M_1 r_1 <1/2$ and let $A:=( \wt 0,
\frac12 r_1)$. We define a function $v$ by
$$
v(x):=\P_x \left(Y_{ \tau^Y_{D \cap B(
0,M_1r_1)}} \in \{y \in D;
y_d
>0\}\right).
$$
By definition $v$ is regular harmonic in $D \cap B(
0,M_1r_1)$ with
respect to $Y$ and vanishes in $D^c$. Applying the function $v$ above to
(\ref{ce0}), we get a Carleson type estimate at $0$, i.e., for any nonnegative
function $u$  which is regular harmonic with respect to $Y$ in $D
\cap B(
0, M_1 r_1)$ and vanishes in $D^c$ we have
\begin{equation}\label{ce1}
u(A) \,\ge \, c^{-1}_{r_1} \frac{v(A)}{v(x)} u(x)  \,\ge \,   c^{-1}_{r_1} v(A)
u(x) \,=\,c_1\,u(x), \quad x \in D \cap B(0, r_1),
\end{equation}
where $c_1=c^{-1}_{r_1} v(A) >0$.
We will construct a bounded positive function $u$ on $\R^d$ which is regular harmonic with respect to $Y$
in $D \cap B(
0, M_1 r)$  and vanishes
in $D^c$ for which (\ref{ce1}) fails.

For $n \ge 1$, we put
\begin{eqnarray*}
C_n&:=&\left\{ (
\wt x,
x_d) \in D: \quad |
\wt x| \le 2^{-n-3} r_1,
\quad  x_2
\le -1+2^{-n}r_1^2 \right\},\\
D_n&:=&\left\{ (
\wt y,
y_d) \in D: \quad
y_d>0,  \quad
|x-y| <1 \quad \mbox{ for some } x \in  C_n\right\}.
\end{eqnarray*}
 It is easy to see that
\begin{equation}\label{ce2}
\overline{D_n} \subset \{(
\wt y,
y_d): |
\wt y| \le ( 2^{-n-3}+ 2^{-(n-1)/2}) r_1, 0 \le
y_d \le 2^{-n}r_1^2\} \subset B(
0,r_1) \cap D,  \quad \mbox{ for } n \ge 2.
\end{equation}
Indeed,  for any $y\in \overline{D_n}$, we have $
y_d\in [0, 2^{-n}r_1^2]$ and
$|y-x |\le 1$ for some $x \in C_n$.
 If  $|
\wt y| > ( 2^{-n-3}+ 2^{-(n-1)/2}) r_1$, $
y_d \ge 0$ and $x \in C_n$, then
$$
|x-y|^2 \ge
x_d^2 + (|
\wt y|-|
\wt x|)^2 \ge (1-2^{-n}r_1^2)^2+ 2^{-(n-1)}r_1^2 >1.$$
Thus, in this case $y \notin  \overline{D_n}$.

For any $n$, let $T^Y_{D_n}$ be the first hitting time of $D_n$ by the
process $Y$. By \eqref{ce2}
$$
\P_A\left(\tau^Y_{D \cap B(
0,M_1r_1)} > T^Y_{D_n}\right) \,\to\,
\P_A\left(\tau^Y_{D \cap B(
0,M_1r_1)} > T^Y_{\{
0\}}\right)
=0, \quad
\mbox{ as } n \to \infty.
$$
Fix $n_0\ge 2$ large so that
\begin{equation}\label{ce3}
\P_A\left(\tau^Y_{D \cap B(
0,M_1r_1)} > T^Y_{D_{n_0}}\right) \, < \, \frac{c_1}{2}
 \end{equation}
and define
$$
u(x)\,:=\P_x \left( Y_{ \tau^Y_{D \cap B(
0,M_1r_1)}} \in C_{n_0}\right).
$$
$u$ is a nonnegative bounded function which is
regular harmonic in $D \cap B(
0,M_1r_1)$ with respect to
$Y$ and vanishes in $D^c$. It also vanishes
continuously on $\partial D \cap B(
0,M_1r_1)$.
Note that by \eqref{Poisson},
$$\P_{A} \left( Y_{ \tau^Y_{D \cap B(
0,M_1r_1)}} \in C_{n_0},
\, \tau^Y_{D \cap B(
0,M_1r_1)} \le T^Y_{D_{n_0}}   \right)
\, =\,
\P_{A} \left(    Y_{ \tau^Y_{D \cap B(
0,M_1r_1)
\setminus D_{n_0}}} \in C_{n_0}      \right)=0.
$$
Thus by the strong Markov property,
\begin{eqnarray*}
u(A) &=&\P_A \left( Y_{ \tau^Y_{D \cap B(
0,M_1r_1)}} \in C_{n_0},\,
\tau^Y_{D \cap B(
0,M_1r_1)} > T^Y_{D_{n_0}} \right)\\
&=&\E_A \left[\P_{Y_{T^Y_{D_{n_0}}}}
\left( Y_{ \tau^Y_{D \cap B(
0,M_1r_1)}} \in C_{n_0} \right);\,
\tau^Y_{D \cap B(
0,M_1r_1)} > T^Y_{D_{n_0}}   \right]\\
&\le&\P_A\left(\tau^Y_{D \cap B(
0,M_1r_1)} > T^Y_{D_{n_0}}\right) \left(
\sup_{x \in D_{n_0}} u(x)\right) \, < \, \frac{c_1}2 \left(\sup_{x \in
D \cap B(
0,r_1)   } u(x)\right).
\end{eqnarray*}
In the last inequality above, we have used (\ref{ce2})--(\ref{ce3}).
But by (\ref{ce1}), $u(A) \,\ge \,c_1   \,  \sup_{x \in D \cap
B(
0,r_1)   }  u(x)$, which gives a contradiction. Thus the boundary
Harnack principle is not true for $D$ at the origin.

By smoothing off the corners of $D$, we can easily construct a bounded
$C^{1, 1}$ domain on which the boundary Harnack principle for $Y$ fails at 0.

\bigskip

\section{Proofs of Theorems \ref{t-main-green} and \ref{t-main-martin}}
As already said in the introduction, once the boundary Harnack principle has been established,
the proofs of Theorems \ref{t-main-green} and \ref{t-main-martin} are similar to the
corresponding proofs in \cite{CKSV2} for the operator $\Delta+ a^{\alpha} \Delta^{\alpha/2}$.
In fact, the proof are even simpler, because \cite{CKSV2} strives for uniformity in the weight $a$.

The proof of Theorem \ref{t-main-green} in the case $d\geq 3$ is by now quite standard.
Once the interior estimates are
established, the full estimates in connected $C^{1,1}$ open sets follow from
the boundary Harnack principle by the method developed
by Bogdan \cite{Bo1} and Hansen \cite{H}. For the operator $\Delta+ a \Delta^{\alpha/2}$
this is accomplished in \cite[Section 3]{CKSV2}. In the present setting the proof from
\cite{CKSV2} carries over almost verbatim. In several places in \cite{CKSV2} one refers
to the form of the L\'evy density, but in fact,
the form of the L\'evy density is only used
to establish uniformity in the weight $a$.

When $d=2$, the above method ceases to work due to the nature of the logarithmic
potential associated with the Laplacian. The proof in \cite[Section 4]{CKSV2} for
the operator $\Delta+ a \Delta^{\alpha/2}$ uses a capacitary argument to derive
the interior upper bound estimate for the Green function. By a scaling consideration
and applying the boundary Harnack principle, one gets sharp Green function upper bound
estimates. For the lower bound estimates, \cite{CKSV2} compares the process with the
subordinate killed Brownian motion when
$D$ is connected, and then extend it to general bounded $C^{1,1}$ by using the
jumping structure of the process. In the present setting, the proof of the lower
bound is exactly the same as in \cite{CKSV2} (see proofs of Theorems 4.2 and 4.4).
The proof of the upper bound is essentially the same as the one in \cite{CKSV2},
except that one has to make several minor modifications. Lemma  4.5 in \cite{CKSV2}
should be replaced by the following statement: There exists $c>0$ such that for any $L>0$,
$$
\mathrm{Cap}^0_{B(0,L)}(\overline{B(0,r)}) \ge
\frac{c}{\log(L/r)}  \quad \text{for every } r \in (0, 3L/4).
$$
This is proved in the same way as \cite[Lemma 4.5]{CKSV2} by using the explicit
formula for the Green function of the ball $B(0,L)\subset \R^2$:
$$
G^0_{B(0,L)}(x,y)=\frac{1}{2\pi}\log\left(1+\frac{
(L^2-|x|^2)(L^2-|y|^2)}{L^2|x-y|^2}\right)\, .
$$
The statement of Lemma 4.6 in \cite{CKSV2} should be changed to: There exists
$c>0$ such that for any $L>0$ and
bounded open set $D$ in $\bR^2$ containing $B(0, L)$
and any $x \in \overline{B(0, \frac{3L}4)}$
$$
G_{D}(x,0)\,\le\, \frac{c}{\mbox{\rm Cap}^0_{D}
\big(\overline{B(0, |x|/2)} \big)} \, \P_x\left(  \sigma_{\overline{B(0, |x|/2 )}} <   \tau_{D}\right) ,
$$
(we refer to \cite{CKSV2} for all unexplained notation).
Next, Corollary 4.7 in \cite{CKSV2} should be replaced by the statement:
There exists $c>0$ such that for any $L>0$ and any $x \in \overline{B(0, 3L/4)}$
$$
G_{B(0,L)}(x,0)\,\le\, c\,\log \left(L/|x| \right).
$$
Finally, the last change is in the proof of Lemma 4.8 in \cite{CKSV2} which
uses a scaling argument. This in our setting can be circumvented by using the
modified statement of \cite[Lemma 4.6]{CKSV2}. The rest of the proof remains
exactly the same.

The proof of Theorem \ref{t-main-martin} is also quite
standard  (see \cite{Bo, CKSV2, KS2, KSV}). In the
current setting we follow step-by-step the proof of the corresponding
result in \cite[Section 6]{CKSV2}.
The main difference is that
\cite{CKSV2} uses the explicit form of the L\'evy density $j^a$ for
the operator $\Delta+ a\Delta^{\alpha/2}$ which is $c(\alpha,d,a)r^{-d-\alpha}$.
This L\'evy density is now replaced by $j$, and it suffice to use properties
\eqref{H:1} and \eqref{H:2} to carry over all arguments. The reader can also compare with
\cite[Section 6]{KSV} where the Martin boundary was identified with the
Euclidean boundary for purely discontinuous processes whose jumping kernel satisfies
\eqref{H:1} and \eqref{H:2}.

\vskip 0.3truein

{\bf Panki Kim}

Department of Mathematical Sciences and Research Institute of Mathematics,
Seoul National University,
Building 27, 1 Gwanak-ro, Gwanak-gu
Seoul 151-747, Republic of Korea

E-mail: \texttt{pkim@snu.ac.kr}

\bigskip

{\bf Renming Song}

Department of Mathematics, University of Illinois, Urbana, IL 61801, USA

E-mail: \texttt{rsong@math.uiuc.edu}

\bigskip

{\bf
Zoran Vondra\v{c}ek}

Department of Mathematics,
University of Zagreb,
Bijeni\v{c}ka c.~30,
Zagreb, Croatia

Email: \texttt{vondra@math.hr}

\end{document}